\documentclass{amsart}
\usepackage{amsthm}
\usepackage[...]{youngtab}
\usepackage{graphicx}
\usepackage{latexsym}
\usepackage{array}
\usepackage{amssymb}
\usepackage{amsmath}
\usepackage{booktabs}
\usepackage{enumitem}
\usepackage{url}
\usepackage{fancyhdr}
\usepackage{txfonts}
\usepackage[all]{xy}
\usepackage{color}
\usepackage{tikz}
\usepackage{caption}
\usepackage{subcaption}
\usepackage{float}
\usetikzlibrary{calc,fadings,decorations.pathreplacing}
\numberwithin{figure}{section}

\newtheorem{thm}{Theorem}[section]
\newtheorem{cor}[thm]{Corollary}
\newtheorem{lem}[thm]{Lemma}

\newtheorem{defn}[thm]{Definition}

\newtheorem{claim}[thm]{Claim}
\newtheorem{notn}[thm]{Notation}

\newcommand\restrict[1]{\raisebox{-.4ex}{$|$}_{#1}}

\begin{document}

\title{Homomorphism Complexes and k-Cores}
\author{Greg Malen}
\address{Department of Mathematics \\
The Ohio State University}
\date{\today}
\maketitle


\begin{abstract}
We prove that the topological connectivity of a graph homomorphism complex Hom($G,K_m$) is at least $m-D(G)-2$, where $\displaystyle D(G)=\max_{H\subseteq G}\delta(H)$. This is a strong generalization of a theorem of \u{C}uki\'{c} and Kozlov, in which $D(G)$ is replaced by the maximum degree $\Delta(G)$. It also generalizes the graph theoretic bound for chromatic number, $\displaystyle\chiup(G)\leq D(G)+1$, as $\displaystyle\chiup(G)=\min\{ m:\text{Hom}(G,K_m)\neq\varnothing\}$. Furthermore, we use this result to examine homological phase transitions in the random polyhedral complexes Hom$(G(n,p),K_m)$ when $p=c/n$ for a fixed constant $c > 0$.
\end{abstract}

\section{Introduction}

The study of graph homomorphism complexes, Hom$(G,H)$ for fixed graphs $G$ and $H$, grew out of a rising interest in finding topological obstructions to graph colorings, as in Lov{\'a}sz's proof of the Kneser Conjecture ~\cite{Lov1}. Specifically, Lov{\'a}sz used the topological connectivity of the Neighborhood Complex, $\mathcal{N}(G)$, to provide a lower bound on the chromatic number of a Kneser graph. It turns out that $\mathcal{N}(G)$ is a specialization of a homomorphism complex, as it is homotopy equivalent to Hom$(K_2,G)$ (see Proposition 4.2 in ~\cite{BabKoz}). The structure of the underlying graphs has many quantifiable effects on the topology of these complexes. The 0-cells of Hom$(G,H)$ are precisely the graph homomorphisms from $G\rightarrow H$, and so the most straighforward such result is that the chromatic number of $G, \ \chiup(G)$, is the minimum $m$ for which Hom$(G,K_m)$ is non-empty.  For $\Delta(G)$ the maximum degree of $G$, Babson and Kozlov conjectured ~\cite{BabKoz}, and \u{C}uki\'{c} and Kozlov later proved ~\cite{CukKoz1} that if $\Delta(G)=d$, then Hom($G,K_m$) is at least $(m-d-2)$-connected. Note that when $m=d+1$ this recovers the well known bound $\chiup(G)\leq\Delta(G)+1$.
\par In graph theory, there is an improvement of this bound using the minimum degree of induced subgraphs, namely that $\displaystyle\chiup(G)\leq\max_{H\subseteq G}\delta(H)+1$, for $\delta(H)$ the minimum degree of $H$. The value $\displaystyle D(G)\coloneqq\max_{H\subseteq G}\delta(H)$ is known as the degeneracy of $G$, a graph property which is commonly utilized in computer science and the study of large networks. Using $D(G)$ in place of $\Delta(G)$, in Section 3 we prove a generalization of this graph theoretic result which will also specialize to the \u{C}uki\'{c}--Kozlov Theorem. In Section 4 we focus on the case that $H=K_3$ and show, with only a few notable exceptions, that Hom($G,K_3$) is disconnected if it is non-empty.
\par Furthermore, in the setting of random graphs the evolution of $D(G)$ for a sparse Erd\H{o}s--R\'{e}nyi graph $G\sim G(n,p)$ has been well studied. For $p=c/n$ for a constant $c > 0$, Pittel, Spencer, and Wormwald exhibited sharp thresholds for the appearance and size of subgraphs with minimum degree $k$ ~\cite{PSWcore}. In Section 5 we combine our results with their work to exhibit both phase transitions and lower bounds for the topological connectivity of the random polyhedral complexes Hom$(G(n,c/n),K_m)$ for constant $c > 0$. When $m=3$, we are able to use the results from Section 4 to evaluate the limiting probability for the connectivity of Hom$(G(n,c/n),K_3)$ for all $c > 0$.
\vskip.5cm
\section{Background and Definitions}

In the following, all graphs are undirected simple graphs, and $K_m$ is the complete graph on $m$ vertices. For a graph $G$ and a vertex $v\in V(G)$, when it is clear from context we may write simply that $v\in G$. In this section we define a graph homomorphism complex, or hom-complex, and various properties which will be used to acquire bounds on the topological connectivity of these complexes.
\vskip.3cm
\begin{defn}
For fixed graphs $G$ and $H$, a graph homomorphism from $G\rightarrow H$ is an edge preserving map on the vertices, $f:V(G)\rightarrow V(H)$ such that $\{ f(v),f(u)\}\in E(H)$ if $\{ v,u\}\in E(G)$. 
\end{defn}

\begin{defn}
For fixed graphs $G$ and $H$, cells in Hom($G,H$) are functions $$\eta: V(G)\rightarrow 2^{V(H)}\setminus \varnothing$$ with the restriction that if $\{ x,y\}\in E(G)$, then $\eta(x)\times\eta(y)\subseteq E(H)$. The dimension of $\eta$ is defined to be $\displaystyle\dim(\eta)=\sum_{v\in V(G)}\left(\left|\eta(v)\right|-1\right)$, with ordering $\eta\subseteq \tau$ if $\eta(v)\subseteq\tau(v)$ for all $v\in V(G)$.
\end{defn}
\vskip.3cm
The set of 0-cells of Hom$(G,H)$ is then precisely the set of graph homomorphisms from $G\rightarrow H$, with higher dimensional cells formed over them as multihomomorphisms. Furthermore, every cell in a hom-complex is a product of simplices, and hom-complexes are entirely determined by their 1-skeletons. If $\eta$ is a product of simplices and its 1-skeleton is contained in Hom$(G,H)$, then $\eta\in\text{Hom}(G,H)$. The 1-skeleton of Hom$(G,H)$ can be thought of in the following manner. Two graph homomorphisms $\eta,\tau:G\rightarrow H$ are adjacent in Hom$(G,H)$ if and only if their images differ on exactly one vertex $v\in G$. So $\eta(v)\neq\tau(v)$, $\eta(u)=\tau(u)$ for all $u\in G\setminus\{v\}$, and the 1-cell joining $\eta$ and $\tau$ is the multihomomorphism $\sigma$ for which $\sigma(v)=\eta(v)\cup\tau(v)$ and $\sigma(u)=\eta(u)=\tau(u)$ for all $u\in G\setminus\{v\}$.
\par Since the property of being a homomorphism is independent on disjoint connected components of a graph, hom-complexes obey a product rule for disjoint unions. For any graphs $G,G'$ and $H$, $$\text{Hom}(G \ {\big. \amalg} \ G',H)=\text{Hom}(G,H)\times\text{Hom}(G',H)$$ Thus, when convenient we can always restrict our attention to connected graphs. Here we also introduce the notion of a graph folding.
\vskip.3cm
\begin{defn}
Denote the neighborhood of a vertex $v\in G$ by $N(v)=\{w\in G:w\sim v\}$. Let $v,u\in G$ be distinct vertices such that $N(v)\subseteq N(u)$. Then a fold of $G$ is a homomorphism from $G\rightarrow G\setminus\{ v\}$ which sends $v\mapsto u$ and is otherwise the identity.
\end{defn}

\begin{lem}
(Babson, Kozlov; Proposition 5.1 in ~\cite{BabKoz}). If $G$ and $H$ are graphs and $v,u\in G$ are distinct vertices such that $N(v)\subseteq N(u)$, then the folding $G\rightarrow G\setminus\{ v\}$ which sends $v\mapsto u$ induces a homotopy equivalence Hom$(G\setminus\{v\},H)\rightarrow\text{Hom}(G,H)$. 
\end{lem}
\vskip.3cm
If $T$ is a tree, for example, then $T$ folds to a single edge, so Hom$(T,K_m)\simeq$ Hom$(K_2,K_m)$. And by Proposition 4.5 in ~\cite{BabKoz}, Hom$(K_n,K_m)$ is homotopy equivalent to a wedge of $(m-n)$-spheres, and in particular Hom$(K_2,K_m)\simeq S^{m-2}$. So Hom$(T,K_m)\simeq S^{m-2}$ for any tree $T$. For a thorough introduction to hom-complexes, see ~\cite{BabKoz}.\\

\begin{notn}
For a graph $G$ define the maximum degree to be $\displaystyle\Delta(G)\coloneqq\max_{v\in V(G)}\{\deg(v)\}$, and the minimum degree to be $\displaystyle\delta(G)\coloneqq\min_{v\in V(G)}\{\deg(v)\}$.
\end{notn}

\begin{defn}
The $k$-core of a graph $G$ is the subgraph $c_k(G)\subseteq G$ obtained by the process of deleting vertices with degree less then k, along with all incident edges, one at a time until there are no vertices with degree less than $k$.
\end{defn}
\vskip.2cm
Regardless of the order in which vertices are deleted, this process always terminates in the unique induced subgraph $c_k(G)\subseteq G$ which is maximal over all subgraphs of $G$ which have minimum degree at least $k$. The $k$-core of a graph may be the empty graph, and the existence of a non-empty $k$-core is a monotone question, as $c_l(G)\subseteq c_j(G)$ whenever $j\leq l$. Here we are interested in the most highly connected non-trivial subgraph.
\vskip.3cm
\begin{defn}
The degeneracy of $G$ is $\displaystyle D(G)\coloneqq\max_{H\subseteq G}\delta(H)$, for $H$ an induced subgraph. 
\end{defn}
\vskip.2cm
Then $D(G)$ is the maximum $k$ such that $c_k(G)$ is non-empty. Given these definitions, we are now able to state the main theorem:
\vskip.3cm
\begin{thm}
For any graph $G$, Hom$(G,K_m)$ is at least $(m-D(G)-2)$-connected.
\end{thm}
\vskip.2cm
And $D(G)\leq\Delta(G)$, hence this will imply the \u{C}uki\'{c}--Kozlov Theorem. Our interest is primarily in applying this result to the case that $G$ is a sparse Erd\H{o}s--R\'{e}nyi random graph, where $D(G)$ has been studied extensively and is much smaller than $\Delta(G)$. It should also be noted that Engstr{\"o}m gave a similar strengthening of the \u{C}uki\'{c}--Kozlov Theorem in ~\cite{Eng1}, replacing $\Delta(G)$ with a graph property which is independent of $D(G)$ and which may provide a better tool for studying the dense regime.\\
\vskip.7cm
\section{Proof of Theorem 2.8}

To prove the theorem, we first introduce a result of Csorba which finds subcomplexes that are homotopy equivalent to Hom$(G,K_m)$ by removing independent sets of $G$.
\vskip 0.3cm
\begin{notn}
Let $G$ and $H$ be graphs. Define $$\text{Ind}(G)\coloneqq\{ S\subseteq V(G) : S \text{ is an independent set}\}$$ And for $I\in \text{Ind}(G)$ define $$\text{Hom}_I(G,H)\coloneqq\left\{\eta\in \text{Hom}(G\setminus I,H): \text{ there is } \overline{\eta}\in \text{Hom}(G,H) \text{ with } \overline{\eta}\restrict{G\setminus I}=\eta\right\}$$
\end{notn}
\vskip0.2cm
So Hom$_I(G,H)$ is the subcomplex of Hom$(G\setminus I,H)$ comprised of all multihomomorphisms from $G\setminus I\rightarrow H$ which extend to multihomomorphisms from $G\rightarrow H$.
\vskip .3cm
\begin{thm}
(Csorba; Theorem 2.36 in ~\cite{Csorba1}) For any graphs $G$ and $H$, and any $I\in\text{Ind}(G)$, Hom$(G,H)$ is homotopy equivalent to Hom$_I(G,H)$.
\end{thm}
\vskip 0.2cm
The proof of this is an application of both the Nerve Lemma and the Quillen Fiber Lemma, and a more thorough examination of this property is given by Schultz in ~\cite{Schultz1}. Note that when $H=K_m$ and $I=\{ v\}$ for any vertex $v\in G$, a multihomomorphism $\eta:G\setminus\{v\}\rightarrow K_m$ has an extension $\overline{\eta}:G\rightarrow K_m$ as long as there is some vertex in $K_m$ which is not in $\eta(w)$ for any $w\in N(v)$. Thus  $$\text{Hom}_{\{ v\}}(G,K_m) = \left\{ \eta\in \text{ Hom}(G\setminus\{ v\},K_m) : \left|\bigcup_{w\in N(v)}\eta(w)\right|\leq m-1\right\}$$.

\begin{lem}
Let G be a graph with a vertex $v$ such that $\deg(v)\leq k$.  Then the $(m-k-1)$-skeleton of Hom$(G\setminus\{ v\},K_m$) is contained in Hom$_{\{ v\}}(G,K_m)$.
\end{lem}

\noindent{\bf Proof of Lemma 3.3.} Let $\eta\in$ Hom$(G\setminus\{ v\},K_m) \ \setminus$ Hom$_{\{ v\}}(G,K_m)$. Then $$m=\left|\bigcup_{w\in N(v)}\eta(w)\right|\leq\sum_{w\in G\setminus\{ v\}}|\eta(w)|$$ And $$\dim(\eta)=\sum_{w\in G\setminus\{ v\}}\left(|\eta(w)|-1\right)=\left(\sum_{w\in G\setminus\{ v\}}|\eta(w)|\right)-|N(v)| \ \geq \ m-k$$ Hence if $\eta\in\text{Hom}(G\setminus\{ v\},K_m)$ with $\dim(\eta)\leq m-k-1$, then $\eta\in\text{Hom}_{\{ v\}}(G,K_m)$. $_{\square}$
\vskip.3cm
So if Hom($G\setminus\{ v\},K_n)$ is $(m-k-2)$-connected, then Hom$_{\{ v\}}(G,K_n)$ is $(m-k-2)$-connected. Combining this with Csorba's theorem, we have the following corollary and the proof of Theorem 2.8:
\vskip.3cm
\begin{cor}
Let G be a graph with a vertex $v$ such that $\deg(v)\leq k$. If Hom$(G\setminus\{ v\},K_n)$ is $(m-k-2)$-connected, then Hom$(G,K_n)$ is $(m-k-2)$-connected.
\end{cor}
\vskip.2cm
\noindent{\bf Proof of Theorem 2.8.} \ 
Let $G$ be a graph with $D(G)=k$. Then $c_{k+1}(G)$ is the empty graph, and there is a sequence $G=G_0, \ G_i=G_{i-1}\setminus\{ v_i\}$, with $\deg_{G_{i-1}}(v_i)\leq k$, terminating in $G_{|V(G)|}=\varnothing$. So $G_{|V(G)|-1}$ is a single vertex, and Hom$(G_{|V(G)|-1},K_n)=\Delta^{m}$, the $m$-simplex, which is contractible, and thus $(m-k-2)$-connected. Hence, by induction, Hom$(G,K_m)$ is also $(m-k-2)$-connected. $_{\square}$
\vskip.7cm
\section{Hom($G,K_3$)}
When $G$ is a graph such that $\chiup(G)\leq3$, Hom($G,K_3)$ has a particularly nice structure. If $G$ does not have an isolated vertex, then Hom($G,K_3)$ is a cubical complex, and Babson and Kozlov showed that it admits a metric with nonpositive curvature ~\cite{BabKoz}. Notice that for a connected graph $G$ with $\chiup(G)\leq3$, the bound obtained by Theorem 2.8 when $m=3$ provides no new information. When $D(G)=1$, $G$ is a tree and folds to a single edge, so Hom$(G,K_3)\simeq\text{Hom}(K_2,K_3)\simeq S^1$. When $D(G)=2$, the bound merely confirms that Hom$(G,K_3)$ is non-empty, and for $D(G) > 2$ it gives no information at all. Here we show that Hom($G,K_3)$ is, in fact, disconnected for a large class of graphs $G$ with $\chiup(G)\leq3$.
\vskip.3cm
\begin{thm}
If $G$ is a graph with $\chiup(G)=3$, then Hom($G,K_3$) is disconnected.
\end{thm}
\vskip0.2cm
This result has been formulated previously in the language of statistical physics, where Glauber dynamics examines precisely the 1-skeleton of Hom$(G(n,p),K_m)$. See, for example, the work of Cereceda, van den Heuvel and Johnson ~\cite{CervdHJ}. We provide a proof here in the context of work on hom-complexes of cycles done by \u{C}uki\'{c} and Kozlov, and offer a more general approach for lifting disconnected components via subgraphs in the following lemma.
\vskip0.3cm
\begin{lem}
Let $G \subseteq G'$ and $H$ be graphs such that Hom($G,H$) is non-empty and disconnected. Let $\eta_1,\eta_2\in \text{Hom}(G,H)$ be 0-cells, i.e. graph homomorphisms from $G\rightarrow H$, such that they are in distinct connected components of Hom$(G,H)$. If there are extensions $\overline{\eta_1},\overline{\eta_2}\in$ Hom$(G',H)$ with $\overline{\eta_i}\restrict{G}=\eta_i$ \ for $i\in\{1,2\}$, then Hom$(G',H)$ is also disconnected.
\end{lem}
\vskip0.2cm
\noindent{\bf Proof of Lemma 4.2. } Let $\eta_1,\eta_2$ be graph homomorphisms from $G\rightarrow H$ such that they are in distinct connected components of Hom$(G,H)$, with extensions $\overline{\eta_1},\overline{\eta_2}\in$ Hom$(G',H)$. We may assume that $\dim(\overline{\eta_i})=0$ for each $i$, since otherwise any 0-cells they contain would also be extensions of $\eta_i$. Suppose that $\overline{\eta_1},\overline{\eta_2}$ are in the same connected component of Hom$(G',H)$. Then there is a path $\overline{\eta_1}=\overline{\tau_0}\sim\overline{\tau_1}\sim\overline{\tau_2}\sim\ldots\sim\overline{\tau_l}=\overline{\eta_2}$ in Hom($G',H$), with $\dim(\overline{\tau_j})=0$ for all $j$. Let $\tau_j = \overline{\tau_j}\restrict{G}$, for $0\leq j\leq l$. Then as functions on $V(G')$, for each $j$ there is one vertex $v_j\in G'$ for which $\overline{\tau_j}(v_j)\neq\overline{\tau_{j+1}}(v_j)$, and they agree on all other vertices. If $v_j\in G'\setminus G$, then for the restrictions $\tau_j=\tau_{j+1}$. If $v_j\in G$, then $\tau_j(v_j)\neq\tau_{j+1}(v_j)$, but they agree on all other vertices, so $\tau_j\sim\tau_{j+1}$ in Hom($G,H$). Hence the path in Hom($G',H$) projects onto a possibly shorter path from $\eta_1$ to $\eta_2$ in Hom($G,H$), which is a contradiction.\\
\\
Therefore the extensions $\overline{\eta_1}$ and $\overline{\eta_2}$ must be in different connected components of Hom($G',H$), which is thus disconnected. $_{\square}$
\vskip0.3cm
The subgraphs we examine to obtain disconnected components of Hom($G,K_3$) will be cycles. In ~\cite{CukKoz2}, \u{C}uki\'{c} and Kozlov fully characterized the homotopy type of Hom($C_n,C_m$) for all $n,m\in\mathbb{N}$. In particular, they showed that all 0-cells in a given connected component have the same number of what they call return points. For these complexes, let $V(C_n)=\{ v_1,\ldots, v_n\}$ such that $v_i\sim v_{(i+1) \text{ mod }n}$ for all $1\leq i\leq n$, and let $V(C_m)=\{ 1,\ldots, m\}$ such that $j\sim (j+1) \text{ mod } m$ for all $1\leq j\leq m$. For the purpose of defining the return number of a 0-cell $\eta$, we momentarily drop the set bracket notation and write $\eta(v_i)=j$. 
\vskip0.3cm
\begin{defn}
A return point of a 0-cell $\eta\in$ Hom($C_n, C_m$) is a vertex $v_i\in C_n$ such that $\eta(v_{i+1})-\eta(v_i)\equiv -1 \text{ mod } m$. Then $r(\eta)$, the return number of $\eta$, is the number of $v_i\in C_n$ which are return points of $\eta$.
\end{defn}
\vskip.2cm
Note that since $\eta$ is a homomorphism, the quantity $\eta(v_{i+1})-\eta(v_i)$ is always $\pm1$ mod $m$. It simply measures in which direction $C_n$ is wrapping around $C_m$ on the edge $\{ v_i,v_{i+1}\}$.
\vskip.3cm
\begin{lem}
(\u{C}uki\'{c}, Kozlov, Lemma 5.3 in ~\cite{CukKoz2}) If two 0-cells of Hom$(C_n,C_m)$ are in the same connected component, then they have the same return number.
\end{lem}
\vskip.2cm
For a 0-cell $\eta\in$ Hom($G,K_m$), one may always obtain another 0-cell by swapping the inverse images of two vertices in $K_m$. When the target graph is a complete graph, we refer to the inverse image of a vertex in $K_m$ as a color class of $\eta$. When $m=3$, we have $K_3=C_3$ and we can track the effect that interchanging two color classes has on the return number. 
\vskip.3cm
\begin{notn}
For a 0-cell $\eta\in$ Hom$(G,K_m)$, let $\eta_{\{ l,j\}}$ denote the $(l,j)$-interchanging 0-cell obtained by defining $\eta^{-1}_{\{ l,j\}}(l)=\eta^{-1}(j), \ \eta^{-1}_{\{ l,j\}}(j)=\eta^{-1}(l)$, and $\eta^{-1}_{\{ l,j\}}(i)=\eta^{-1}(i)$ for all $i\notin\{ l,j\}$.
\end{notn}
\begin{lem}
Fix a pair $\{ l,j\}\subset\{1,2,3\}, \ l\neq j$. Then for any 0-cell $\eta\in$ Hom($C_n,K_3$), every $v_i\in C_n$ is a return point for exactly one of $\eta$ and $\eta_{\{ l,j\}}$. Hence $r\left(\eta_{\{ l,j\}}\right)=n-r(\eta)$.
\end{lem}
\vskip.2cm
\noindent {\bf Proof of Lemma 4.6.} Consider the edge $\{ v_i,v_{i+1}\}$. SInce the target graph is $K_3$, we have that $\{\eta(v_i),\eta(v_{i+1})\}\cap\{ l,j\}\neq\varnothing$. If $\{\eta(v_i),\eta(v_{i+1})\}=\{ l,j\}$, then $$\eta_{\{ l,j\}}(v_{i+1})-\eta_{\{ l,j\}}(v_i)=-(\eta(v_{i+1})-\eta(v_i))$$ Thus $v_i$ is a return point of $\eta$ if and only if it is not a return point of $\eta_{\{ l,j\}}$. Alternatively, if $|\{\eta(v_i),\eta_(v_{i+1}\}\cap\{ l,j\}|=1$, then one of $v_i$ and $v_{i+1}$ has a stationary image under an $(l,j)$-interchange while the other does not. By fixing the image of one vertex and changing the other, this flips the direction that $C_n$ is wrapping around $K_3$ on the edge $\{ v_i,v_{i+1}\}$. So returns become non-returns and vice versa. Therefore each $v_i\in C_n$ is a return point of exactly one of $\eta$ and $\eta_{\{ l,j\}}$ for any fixed $\{ l,j\}\subset\{1,2,3\}$. $_{\square}$
\vskip.3cm
The proof of Theorem 4.1 will then proceed by finding an odd cycle in $G$ and using a color class interchange to produce 0-cells in distinct connected components, which can be lifted from Hom($C_{2k+1},K_3$) to Hom($G,K_3$).
\vskip.4cm
\noindent {\bf Proof of Theorem 4.1} Since $\chiup(G)=3$, Hom($G,K_3$) is non-empty and $G$ contains an odd cycle $H=C_{2k+1}$ for some $k\in\mathbb{N}$. Note that we do not require $H$ to be an induced cycle. Let $\eta\in$ Hom($G,K_3$) be a 0-cell, and let $\tau=\eta\restrict{H}$ be the induced homomorphism on $H$. Fix a distinct pair $\{ l,j\}\subset \{1,2,3\}$ and consider the $(l,j)$-interchange $\eta_{\{ l,j\}}$. It is straightforward that $\eta_{\{ l,j\}}\restrict{H}=\tau_{\{ l,j\}}\in$ Hom$(H,K_3)$. Thus by Lemma 4.6, $r\left(\tau_{\{ l,j\}}\right)=2k+1-r(\tau)\neq r(\tau)$. Hence, by lemma 4.4 they are in different connected components of Hom$(H,K_3)$, and by Lemma 4.2 their extensions $\eta$ and $\eta_{\{ l,j\}}$ are in different connected components of Hom($G,K_3$), which is thus disconnected. $_{\square}$
\vskip.4cm
Bipartite graphs are a bit more complicated. We have seen that for any graph $G$ which folds to a single edge, Hom$(G,K_3)\simeq S^1$. When $G$ is bipartite and does not fold to an edge, it must contain an even cycle $C_{2k}$ for $k\geq3$. So by the same method one can start with $\eta\in$ Hom$\left(G,K_3\right)$ and take its restriction $\eta\restrict{C_{2k}}=\tau\in$ Hom$(C_{2k},K_3$). But if $r(\tau)=k$, then interchanging two color classes no longer guarantees disjoint components. In particular, the existence of too many 4-cycles close to $C_{2k}$ may force $r(\tau)=k$. For example, for $Q_3$ the 1-skeleton of the 3-cube, Hom($Q_3,K_3$) is connected. Note that $D(Q_3)=3$, so the lower bound achieved by Theorem 2.8 does not provide any information. Restricting to the case that $G$ omits the subgraphs in Figure 4.1 is sufficient to ensure that this does not happen, and that Hom$(G,K_3)$ is disconnected.

\begin{figure}[H]
\hspace*{2.6cm}
   \begin{subfigure}[t]{\textwidth}
        \begin{tikzpicture}
        \coordinate (A1) at (0cm,-.5cm);
	\coordinate (A2) at (0cm,.5cm);
	\coordinate (A3) at (-1cm,-.5cm);
	\coordinate (A4) at (-1cm,.5cm);
	\coordinate (A5) at (1cm,-.5cm);
	\coordinate (A6) at (1cm,.5cm);
	\draw[thick] (A1) -- (A2);
	\draw[thick] (A1) -- (A3);
	\draw[thick] (A3) -- (A4)
		node at (barycentric cs:A3=1,A4=1) [left] {$H_1 =$};
	\draw[thick] (A1) -- (A5);
	\draw[thick] (A5) -- (A6);
	\draw[thick] (A6) -- (A2);
	\draw[thick] (A2) -- (A4);
	\draw[fill=black] (A1) circle (0.15em);
	\draw[fill=black] (A2) circle (0.15em);
	\draw[fill=black] (A3) circle (0.15em);
	\draw[fill=black] (A4) circle (0.15em);
	\draw[fill=black] (A5) circle (0.15em);
	\draw[fill=black] (A6) circle (0.15em);

        \coordinate (B1) at (4.2cm,0cm);
	\coordinate (B2) at (4.2cm,.5cm);
	\coordinate (B3) at (4.2cm,-.5cm);
	\coordinate (B4) at (3.4cm,0cm);
	\coordinate (B5) at (5cm,0cm);
	\draw[thick] (B1) -- (B2);
	\draw[thick] (B1) -- (B3);
	\draw[thick] (B3) -- (B4);
	\draw[thick] (B2) -- (B4);
	\draw[thick] (B3) -- (B5);
	\draw[thick] (B2) -- (B5);
	\draw[fill=black] (B1) circle (0.15em);
	\draw[fill=black] (B2) circle (0.15em);
	\draw[fill=black] (B3) circle (0.15em);
	\draw[fill=black] (B4) circle (0.15em)
		node[left] {$H_2 = \ $};
	\draw[fill=black] (B5) circle (0.15em);

\end{tikzpicture}
\end{subfigure}
\caption{If Hom$(G,K_3)$ is connected, then $G$ must contain $H_1$ or $H_2$.}
\end{figure}
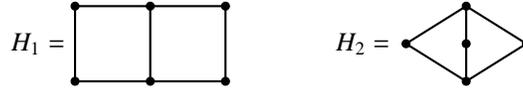
	
\begin{thm}
If $G$ is a finite, bipartite, connected graph which does not fold to an edge and does not contain $H_1$ or $H_2$ as a subgraph, then Hom($G,K_3$) is disconnected.
\end{thm}
\vskip0.2cm
To show this requires a great deal of case-by-case structural analysis when a minimal even cycle in $G$ has length 6, 8, or 10. The method of the proof will suggest a stronger, albeit more technical result.\\
\vskip.2cm
\noindent {\bf Proof of Theorem 4.7.} Let $G$ be a finite, bipartite, connected graph which does not fold to an edge and does not contain $H_1$ or $H_2$ as a subgraph. Since $G$ is finite, bipartite and does not fold to an edge, it must contain an even cycle of length greater than 4. Let $H=C_{2k}, \ k\geq3$ be a fixed cycle in $G$ which has minimal length over all cycles in $G$ of length greater than 4. Label the vertices of $H$ as $v_1,\ldots,v_{2k}$ such that $v_i\sim v_{(i+1)\text{ mod } 2k}$ for all $i$. $H$ cannot have internal chords, as any such edge would create a cycle $C_l$ for $4 < l < 2k$, or in the case that $k=3$ an antipodal chord would create a copy of $H_1$. For $2k\equiv i$ mod 3, take $\tau\in$ Hom($H,K_3)$ such that $\tau(v_j)=\{ j \text{ mod } 3\}$ for $j\leq2k-i$, with $\tau(v_{2k})=\{2\}$ if $i=1$, and $\tau(v_{2k-1})=\{1\}, \ \tau(v_{2k})=\{2\}$ if $i=2$. Figure 4.2 depicts $\tau$ for $k\in\{3,4,5\}$.
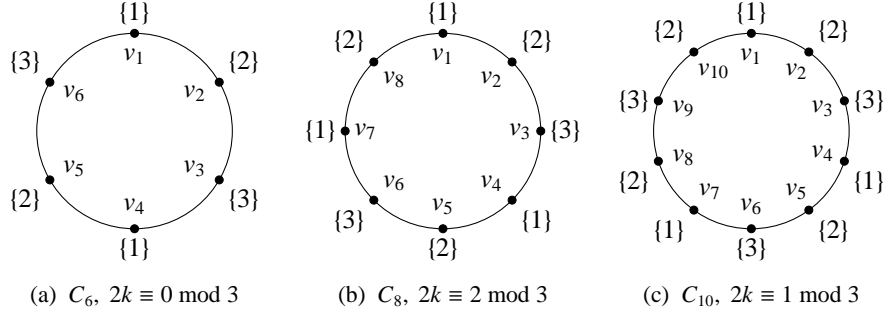
\begin{figure}[H]
\hspace*{.4cm}
\begin{subfigure}[t]{\textwidth}
\begin{tikzpicture}
\draw (-4.1,0) circle (1.3cm) node at (-4.1,-2.2) {\small (a) \ $C_6, \ 2k\equiv 0\text{ mod }3$};
\coordinate (O1) at (-4.1,0);
\coordinate (a1) at (-4.1,1.3);
\draw[fill=black] (a1) circle (0.15em) node at (-4.1,1) {$v_1$} node[above] {$\{1\}$};
\coordinate (a2) at ($(O1)!1!-60:(a1)$);
\draw[fill=black] (a2) circle (0.15em) node at (-3.3,.5) {$v_2$} node[above right] {$\{2\}$};
\coordinate (a3) at ($(O1)!1!-60:(a2)$);
\draw[fill=black] (a3) circle (0.15em) node at (-3.3,-.5) {$v_3$} node[below right] {$\{3\}$};
\coordinate (a4) at ($(O1)!1!-60:(a3)$);
\draw[fill=black] (a4) circle (0.15em) node at (-4.1,-1) {$v_4$} node[below] {$\{1\}$};
\coordinate (a5) at ($(O1)!1!-60:(a4)$);
\draw[fill=black] (a5) circle (0.15em) node at (-4.9,-.5) {$v_5$} node[below left] {$\{2\}$};
\coordinate (a6) at ($(O1)!1!-60:(a5)$);
\draw[fill=black] (a6) circle (0.15em) node at (-4.9,.5) {$v_6$} node[above left] {$\{3\}$};

\draw (0,0) circle (1.3cm) node at (0,-2.2) {\small (b) \ $C_8, \ 2k\equiv 2\text{ mod }3$};
\coordinate (b3) at (0:1.3);
\draw[fill=black] (b3) circle (0.15em) node[left] {$v_3$} node[right] {$\{3\}$};
\coordinate (b2) at (45:1.3);
\draw[fill=black] (b2) circle (0.15em) node[below left] {$v_2$} node[above right] {$\{2\}$};
\coordinate (b1) at (90:1.3);
\draw[fill=black] (b1) circle (0.15em) node at (0,1) {$v_1$} node[above] {$\{1\}$};
\coordinate (b8) at (135:1.3);
\draw[fill=black] (b8) circle (0.15em) node[below right] {$v_8$} node[above left] {$\{2\}$};
\coordinate (b7) at (180:1.3);
\draw[fill=black] (b7) circle (0.15em) node[right] {$v_7$} node[left] {$\{1\}$};
\coordinate (b6) at (225:1.3);
\draw[fill=black] (b6) circle (0.15em) node[above right] {$v_6$} node[below left] {$\{3\}$};
\coordinate (b5) at (270:1.3);
\draw[fill=black] (b5) circle (0.15em) node at (0,-1) {$v_5$} node[below] {$\{2\}$};
\coordinate (b4) at (315:1.3);
\draw[fill=black] (b4) circle (0.15em) node[above left] {$v_4$} node[below right] {$\{1\}$};

\draw (4.1,0) circle (1.3cm);
\coordinate (O2) at (4.1,0) node at (4.1,-2.2) {\small (c) \ $C_{10}, \ 2k\equiv 1\text{ mod }3$};
\coordinate (v1) at (4.1,1.3);
\draw[fill=black] (v1) circle (0.15em) node at (4.1,1) {$v_1$} node[above] {$\{1\}$};
\coordinate (v2) at ($(O2)!1!-36:(v1)$);
\draw[fill=black] (v2) circle (0.15em) node at (4.7,.8) {$v_2$} node[above right] {$\{2\}$};
\coordinate (v3) at ($(O2)!1!-36:(v2)$);
\draw[fill=black] (v3) circle (0.15em) node at (5.05,.3) {$v_3$} node[right] {$\{3\}$};
\coordinate (v4) at ($(O2)!1!-36:(v3)$);
\draw[fill=black] (v4) circle (0.15em) node at (5.05,-.3) {$v_4$} node[below right] {$\{1\}$};
\coordinate (v5) at ($(O2)!1!-36:(v4)$);
\draw[fill=black] (v5) circle (0.15em) node at (4.7,-.8) {$v_5$} node[below right] {$\{2\}$};
\coordinate (v6) at ($(O2)!1!-36:(v5)$);
\draw[fill=black] (v6) circle (0.15em) node at (4.1,-1) {$v_6$} node[below] {$\{3\}$};
\coordinate (v7) at ($(O2)!1!-36:(v6)$);
\draw[fill=black] (v7) circle (0.15em) node at (3.55,-.8) {$v_7$} node[below left] {$\{1\}$};
\coordinate (v8) at ($(O2)!1!-36:(v7)$);
\draw[fill=black] (v8) circle (0.15em) node at (3.2,-.35) {$v_8$} node[below left] {$\{2\}$};
\coordinate (v9) at ($(O2)!1!-36:(v8)$);
\draw[fill=black] (v9) circle (0.15em) node at (3.2,.3) {$v_9$} node[left] {$\{3\}$};
\coordinate (v10) at ($(O2)!1!-36:(v9)$);
\draw[fill=black] (v10) circle (0.15em) node at (3.6,.8) {$v_{10}$} node[above left] {$\{2\}$};
\end{tikzpicture}
\end{subfigure}
\caption{Image of $\tau$ on $H=C_{2k}$}
\end{figure}
\vskip0.2cm
\noindent As defined, only the last two vertices of $C_{2k}$ can be return points of $\tau$. So $r(\tau)\leq2$ in all cases, and hence is not equal to $k$ for any $k\geq3$. Then for any $\{ l,j\}\subset\{1,2,3\}$, $r(\tau_{\{ l,j\}})\neq r(\tau)$. Thus $\tau_{\{ l,j\}}$ and $\tau$ are in distinct connected components of Hom$(H,K_3)$. So if we can construct an extension $\eta$ to all of $G$, then $\eta_{\{ l,j\}}$ will extend $\tau_{\{ l,j\}}$ and hence $\eta$ and $\eta_{\{ l,j\}}$ will be in distinct connected component of Hom($G,K_3$).\\
\\
Define $d(u,v)$ to be the minimal number of edges in a path connecting two vertices $u,v\in G$. Let $\sigma$ be a bipartition of $G$, viewed as a 0-cell of Hom($G,K_3$) with image contained in $\{1,2\}$. We will let $\eta=\sigma$ on $G\setminus H$ and then make adjustments as necessary when the definition conflicts with $\tau$. Define the following sets:\\
\\
\hspace*{2cm} $B_i \ \ \coloneqq\{ u\in G: \min\{ d(u,v):v\in H\}=i\}$\\
\hspace*{2cm} $A \ \ \ \coloneqq\{ u\in B_1: \exists \ v\in N(u)\cap H$ with $\sigma(u)=\tau(v)\}$\\
\hspace*{2cm} $A_{i,j}\coloneqq\{ u\in A: |N(u)\cap H|=i, |N(u)\cap A|=j\}$\\
\\
So $B_i$ is the set of points distance $i$ from $H$, and $A$ is the set of points in $G\setminus H$ on which we cannot define $\eta=\sigma$. To arrive at an appropriate definition for $\eta$ on the set $A$, we must first make some structural observations concerning the partitioning of $A$ into the subsets $A_{i,j}$.
\vskip.3cm
\begin{claim}
For every $u\in B_1, \ |N(u)\cap H|\leq 2$. And if $|N(u)\cap H|=2$, then $N(u)\cap H=\{ v_i, v_{(i+2)\text{ mod }2k}\}$ for some $i\leq2k$.
\end{claim}
\vskip.2cm
\noindent {\bf Proof of Claim 4.8.} Let $u\in B_1$. Then connecting $u$ to any two points in $H$ creates a cycle of length at most $k+2 < 2k$, as in Figure 4.3 (a). Since $G$ is bipartite and $H$ is minimal with respect to having length at least 6, this is valid only if it creates a 4-cycle, with $N(u)\cap H = \{ v_i, v_{(i+2)\text{ mod }2k}\}$ for some $i\leq2k$. Similarly, any further points of attachment would need to be in a 4-cycle with each of the other two vertices in $N(u)\cap H$. But this is only possible if $H=C_6$, in which case it would create several copies of $H_1$, shown in Figure 4.3 (b). $_{\square}$\\
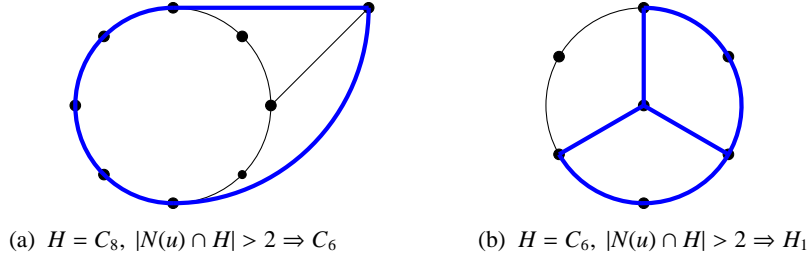
\begin{figure}[H]
\hspace*{.5cm}
\begin{subfigure}[t]{.4\textwidth}
\begin{tikzpicture}
\draw (0,0) circle (1.3cm);
\coordinate (O) at (2.6,1.3);
\draw[fill=black] (O) circle (0.2em);
\coordinate (v3) at (0:1.3);
\draw[fill=black] (v3) circle (0.2em);
\coordinate (v2) at (45:1.3);
\draw[fill=black] (v2) circle (0.2em);
\coordinate (v1) at (90:1.3);
\draw[fill=black] (v1) circle (0.2em);
\coordinate (v8) at (135:1.3);
\draw[fill=black] (v8) circle (0.2em);
\coordinate (v7) at (180:1.3);
\draw[fill=black] (v7) circle (0.2em);
\coordinate (v6) at (225:1.3);
\draw[fill=black] (v6) circle (0.2em);
\coordinate (v5) at (270:1.3);
\draw[fill=black] (v5) circle (0.2em);
\coordinate (v4) at (315:1.3);
\draw[fill=black] (v4) circle (0.15em) node at (0,-1.8) {\small (a) \ $H=C_8, \ |N(u)\cap H| > 2\Rightarrow C_6$};
\draw[blue, ultra thick] (v1) arc (90:270:1.3);
\draw[thin] (O) -- (v3);
\draw[blue, ultra thick] (O) arc (0:-90:2.6);
\draw[blue, ultra thick] (O) -- (v1);
\end{tikzpicture}
\end{subfigure}
\hspace{1cm}
\begin{subfigure}[t]{.4\textwidth}
\begin{tikzpicture}
\draw (0,0) circle (1.3cm);
\coordinate (O) at (0:0);
\draw[fill=black] (O) circle (0.2em);
\coordinate (v1) at (90:1.3);
\draw[fill=black] (v1) circle (0.2em);
\coordinate (v2) at (30:1.3);
\draw[fill=black] (v2) circle (0.2em);
\coordinate (v3) at (330:1.3);
\draw[fill=black] (v3) circle (0.2em);
\coordinate (v4) at (270:1.3);
\draw[fill=black] (v4) circle (0.2em);
\coordinate (v5) at (210:1.3);
\draw[fill=black] (v5) circle (0.2em);
\coordinate (v6) at (150:1.3);
\draw[fill=black] (v6) circle (0.2em) node at (0,-1.8) {\small (b) \ $H=C_6, \ |N(u)\cap H| > 2\Rightarrow H_1$};
\draw[blue, ultra thick] (O) -- (v1);
\draw[blue, ultra thick] (O) -- (v3);
\draw[blue, ultra thick] (O) -- (v5);
\draw[blue, ultra thick] (v5) arc (-150:90:1.3);
\end{tikzpicture}
\end{subfigure}
\caption{Obstructions to $|N(u)\cap H| > 2$}
\end{figure}

\noindent This imposes that $i\leq2$ for $A_{i,j}\subseteq A$, while the following two claims will show that $j$ is at most 1, and that $j=1$ only if $i=1$.
\vskip.3cm
\begin{claim}
Let $k\geq4$, and let $u\in B_1$. Then $|N(u)\cap B_1|\leq1$, and if $N(u)\cap B_1\neq\varnothing$ then $|N(u)\cap H|=1$.
\end{claim}
\vskip.2cm
\noindent {\bf Proof of Claim.} Let $k\geq 4$, let $u\in B_1$, and let $w\in N(u)\cap B_1$, with $u\sim v_i$ and $w\sim v_j$. $G$ is bipartite, so $v_i\neq v_j$. If $v_i\nsim v_j$, as in Figure 4.4 (a), there is a cycle of length $l$ with $4 < l \leq k+3 < 2k$. So $v_i\sim v_j$. Suppose further that $|N(u)\cap H|=2$. Then $N(u)\cap H =\{ v_{(j-1)\text{ mod }2k}, v_{(j+1)\text{ mod }2k}\}$ since both must be adjacent to $v_j$, as in Figure 4.4 (b). But then the induced graph on $\{u,w,v_{(j-1)\text{ mod }2k},v_j,v_{(j+1)\text{ mod }2k}\}$ is a copy of $H_2$. So if $N(u)\cap B_1\neq\varnothing$ for $u\in B_1$, then $|N(u)\cap H| = 1$. Now suppose $|N(u)\cap B_1| > 1$. Then there is a vertex $x\in N(u)\cap B_1$ with $x\neq w$. As with $w$, the vertex in $N(x)\cap H$ must be adjacent to $v_i$. If $x\sim v_j$, as in Figure 4.4 (c), then $\{ u,w,x,v_i,v_j\}$ is a copy of $H_2$, and if $x$ is connected to the other neighbor of $v_i$, as in Figure 4.4 (d), then it creates a copy of $H_1$. Thus $|N(u)\cap B_1| \leq1$. $_\square$
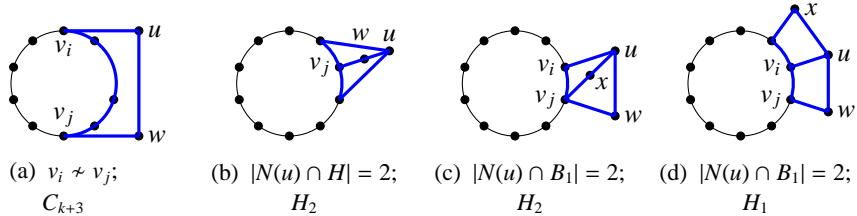
\begin{figure}[H]
\hspace*{.5cm}
\begin{subfigure}[t]{\textwidth}
\begin{tikzpicture}
\draw (-4.5,0) circle (.7cm);
\coordinate (o1) at (-4.5,0);
\coordinate (v1) at (-4.5,.7);
\draw[fill=black] (v1) circle (0.15em) node[below] {$v_i$};
\coordinate (v2) at ($(o1)!1!-36:(v1)$);
\draw[fill=black] (v2) circle (0.15em);
\coordinate (v3) at ($(o1)!1!-72:(v1)$;
\draw[fill=black] (v3) circle (0.15em);
\coordinate (v4) at ($(o1)!1!-108:(v1)$);
\draw[fill=black] (v4) circle (0.15em);
\coordinate (v5) at ($(o1)!1!-144:(v1)$);
\draw[fill=black] (v5) circle (0.15em);
\coordinate (v6) at ($(o1)!1!180:(v1)$);
\draw[fill=black] (v6) circle (0.15em) node[above] {$v_j$};
\coordinate (v7) at ($(o1)!1!144:(v1)$)
	node at (-4.5cm,-1.2cm) {\small (a) \ $v_i\nsim v_j;$};
\draw[fill=black] (v7) circle (0.15em)
	node at (-4.5,-1.6) {\small $ C_{k+3}$};
\coordinate (v8) at ($(o1)!1!108:(v1)$);
\draw[fill=black] (v8) circle (0.15em);
\coordinate (v9) at ($(o1)!1!72:(v1)$:.7);
\draw[fill=black] (v9) circle (0.15em);
\coordinate (v10) at ($(o1)!1!36:(v1)$);
\draw[fill=black] (v10) circle (0.15em);
\coordinate (A) at (-3.5,.7);
\draw[fill=black] (A) circle (0.15em) node[right] {$u$};
\coordinate (B) at (-3.5,-.7);
\draw[fill=black] (B) circle (0.15em) node[right] {$w$};
\draw[blue, very thick] (A) -- (B);
\draw[blue, very thick] (A) -- (v1);
\draw[blue, very thick] (B) -- (v6);
\draw[blue, very thick] (v1) arc (90:-90:.7);

\draw (-1.5,0) circle (.7cm);
\coordinate (o2) at (-1.5,0);
\coordinate (a1) at (-1.5,.7);
\draw[fill=black] (a1) circle (0.15em);
\coordinate (a2) at ($(o2)!1!-36:(a1)$);
\draw[fill=black] (a2) circle (0.15em);
\coordinate (a3) at ($(o2)!1!-72:(a1)$);
\draw[fill=black] (a3) circle (0.15em) node[left] {$v_j$};
\coordinate (a4) at ($(o2)!1!-108:(a1)$);
\draw[fill=black] (a4) circle (0.15em);
\coordinate (a5) at ($(o2)!1!-144:(a1)$);
\draw[fill=black] (a5) circle (0.15em);
\coordinate (a6) at ($(o2)!1!180:(a1)$)
	node at (-1.3cm,-1.2cm) {\small (b) \  $|N(u)\cap H|=2;$};
\draw[fill=black] (a6) circle (0.15em)
	node at (-1.3,-1.6) {\small $H_2$};
\coordinate (a7) at ($(o2)!1!144:(a1)$);
\draw[fill=black] (a7) circle (0.15em);
\coordinate (a8) at ($(o2)!1!108:(a1)$);
\draw[fill=black] (a8) circle (0.15em);
\coordinate (a9) at ($(o2)!1!72:(a1)$);
\draw[fill=black] (a9) circle (0.15em);
\coordinate (a10) at ($(o2)!1!36:(a1)$);
\draw[fill=black] (a10) circle (0.15em);
\coordinate (A) at ($(o2)!1.5!-72:(a1)$);
\draw[fill=black] (A) circle (0.15em) node at (-.55,.65) {$w$};
\coordinate (B) at ($(o2)!2!-72:(a1)$);
\draw[fill=black] (B) circle (0.15em)node[above] {$u$};
\draw[blue, very thick] (A) -- (B);
\draw[blue, very thick] (B) -- (a2);
\draw[blue, very thick] (B) -- (a4);
\draw[blue, very thick] (A) -- (a3);
\draw[blue, very thick] (a2) arc (54:-18:.7);

\draw (1.5,0) circle (.7cm);
\coordinate (o3) at (1.5,0);
\coordinate (b1) at (1.5,.7);
\draw[fill=black] (b1) circle (0.15em);
\coordinate (b2) at ($(o3)!1!-36:(b1)$);
\draw[fill=black] (b2) circle (0.15em);
\coordinate (b3) at ($(o3)!1!-72:(b1)$);
\draw[fill=black] (b3) circle (0.15em) node[left] {$v_i$};
\coordinate (b4) at ($(o3)!1!-108:(b1)$);
\draw[fill=black] (b4) circle (0.15em) node[left] {$v_j$};
\coordinate (b5) at ($(o3)!1!-144:(b1)$);
\draw[fill=black] (b5) circle (0.15em);
\coordinate (b6) at ($(o3)!1!180:(b1)$)
	node at (1.7cm,-1.2cm) {\small (c) \ $|N(u)\cap B_1|=2;$};
\draw[fill=black] (b6) circle (0.15em)
	node at (1.7cm,-1.6cm) {\small $H_2$};
\coordinate (b7) at ($(o3)!1!144:(b1)$);
\draw[fill=black] (b7) circle (0.15em);
\coordinate (b8) at ($(o3)!1!108:(b1)$);
\draw[fill=black] (b8) circle (0.15em);
\coordinate (b9) at ($(o3)!1!72:(b1)$);
\draw[fill=black] (b9) circle (0.15em);
\coordinate (b10) at ($(o3)!1!36:(b1)$);
\draw[fill=black] (b10) circle (0.15em);
\coordinate (A) at ($(o3)!2!-72:(b1)$);
\draw[fill=black] (A) circle (0.15em) node[right] {$u$};
\coordinate (B) at ($(o3)!2!-108:(b1)$);
\draw[fill=black] (B) circle (0.15em) node[right] {$w$};
\coordinate (C) at ($(b4)!.5!(A)$);
\draw[fill=black] (C) circle (0.15em) node at (2.65,0) {$x$};
\draw[blue, very thick] (A) -- (B);
\draw[blue, very thick] (A) -- (C);
\draw[blue, very thick] (C) -- (b4);
\draw[blue, very thick] (A) -- (b3);
\draw[blue, very thick] (B) -- (b4);
\draw[blue, very thick] (b3) arc (18:-18:.7);

\draw (4.5,0) circle (.7cm);
\coordinate (o4) at (4.5,0); 
\coordinate (c1) at (4.5,.7);
\draw[fill=black] (c1) circle (0.15em);
\coordinate (c2) at ($(o4)!1!-36:(c1)$);
\draw[fill=black] (c2) circle (0.15em);
\coordinate (c3) at ($(o4)!1!-72:(c1)$);
\draw[fill=black] (c3) circle (0.15em) node[left] {$v_i$};
\coordinate (c4) at ($(o4)!1!-108:(c1)$);
\draw[fill=black] (c4) circle (0.15em) node[left] {$v_j$};
\coordinate (c5) at ($(o4)!1!-144:(c1)$);
\draw[fill=black] (c5) circle (0.15em);
\coordinate (c6) at ($(o4)!1!180:(c1)$)
	node at (4.7cm,-1.2cm) {\small (d) \  $|N(u)\cap B_1|=2;$};
\draw[fill=black] (c6) circle (0.15em)
	node at (4.7,-1.6) {\small $H_1$};
\coordinate (c7) at ($(o4)!1!144:(c1)$);
\draw[fill=black] (c7) circle (0.15em);
\coordinate (c8) at ($(o4)!1!108:(c1)$);
\draw[fill=black] (c8) circle (0.15em);
\coordinate (c9) at ($(o4)!1!72:(c1)$);
\draw[fill=black] (c9) circle (0.15em);
\coordinate (c10) at ($(o4)!1!36:(c1)$);
\draw[fill=black] (c10) circle (0.15em);
\coordinate (A) at ($(o4)!1.75!-72:(c1)$);
\draw[fill=black] (A) circle (0.15em) node[right] {$u$};
\coordinate (B) at ($(o4)!1.75!-108:(c1)$);
\draw[fill=black] (B) circle (0.15em) node[right] {$w$};
\coordinate (C) at ($(o4)!1.75!-36:(c1)$);
\draw[fill=black] (C) circle (0.15em) node[right] {$x$};
\draw[blue, very thick] (A) -- (B);
\draw[blue, very thick] (A) -- (C);
\draw[blue, very thick] (C) -- (c2);
\draw[blue, very thick] (A) -- (c3);
\draw[blue, very thick] (B) -- (c4);
\draw[blue, very thick] (c2) arc (54:-18:.7);
\end{tikzpicture}
\end{subfigure}
\caption{Invalid structures in $G$ if Claim 4.9 is false.}
\end{figure}

\noindent For $k=3$ a similar result holds when $B_1$ is replaced by the subset $A$:
\vskip.2cm
\begin{claim}
Let $k=3$, and let $u\in A$. Then $|N(u)\cap A|\leq1$, and if $N(u)\cap A\neq\varnothing$ then $|N(u)\cap H|=1$.
\end{claim}
\vskip.2cm
\noindent {\bf Proof of Claim.} Let $u\in A$ and let $w\in N(u)\cap A$. Recall that for $2k=6, \ \tau(v_i)=\{ i\text{ mod } 3\}$ for all $i$. Then without loss of generality assume that $\sigma(u)=\{1\}$ and $u\sim v_1$. So $\sigma(w)=\{2\}$, and since $w\in A$ it must be adjacent to at least one of $v_2$ and $v_5$. But $G$ is bipartite and $u\sim v_5$ creates a 5-cycle, as indicated in Figure 4.5 (a). Thus $w\sim v_2$, and the same argument holds to show that $u\sim v_1$ if one first assumes that $w\sim v_2$. Furthermore, if $w\sim v_4$ or $w\sim v_6$ there is a copy of $H_1$ or $H_2$, respectively, depicted in Figure 4.5 (b), (c). By symmetry $u\nsim v_3, \ u\nsim v_5$. Hence, if $u\in A$ with $N(u)\cap A\neq\varnothing$, then $|N(u)\cap H|=1$. Now suppose there is $x\in N(u)\cap A$ with $x\neq w$. Then $\sigma(x)=\{2\}$, and by the preceding argument we have that $x\sim v_2$. But, as shown in Figure 4.5 (d), $\{ u,w,x,v_1,v_2\}$ is then a copy of $H_2$. Hence for $u\in A, \ |N(u)\cap A|\leq1$. $_{\square}$
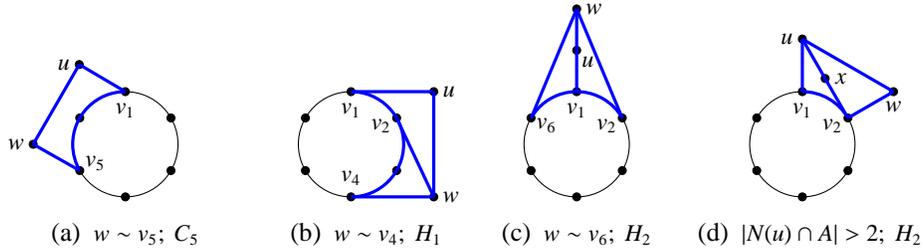
\begin{figure}[H]
\begin{subfigure}[t]{\textwidth}
\begin{tikzpicture}
\draw (-4.5,0) circle (.7cm)
	node at (-4.5,-1.2) {(a) \ \small $w\sim v_5; \ C_5$};
\coordinate (o1) at (-4.5,0);
\coordinate (v1) at (-4.5,.7);
\draw[fill=black] (v1) circle (0.15em) node at (-4.5cm,.45cm) {\small $v_1$};
\coordinate (v2) at ($(o1)!1!-60:(v1)$);
\draw[fill=black] (v2) circle (0.15em);
\coordinate (v3) at ($(o1)!1!-120:(v1)$);
\draw[fill=black] (v3) circle (0.15em);
\coordinate (v4) at ($(o1)!1!180:(v1)$);
\draw[fill=black] (v4) circle (0.15em);
\coordinate (v5) at ($(o1)!1!120:(v1)$);
\draw[fill=black] (v5) circle (0.15em) node at (-4.9cm,-.25cm) {\small $v_5$};
\coordinate (v6) at ($(o1)!1!60:(v1)$);
\draw[fill=black] (v6) circle (0.15em);
\coordinate (A) at ($(o1)!1.75!30:(v1)$);
\draw[fill=black] (A) circle (0.15em) node[left] {\small $u$};
\coordinate (B) at ($(o1)!1.75!90:(v1)$);
\draw[fill=black] (B) circle (0.15em) node[left] {\small $w$};
\draw[blue, very thick] (A) -- (B);
\draw[blue, very thick] (A) -- (v1);
\draw[blue, very thick] (B) -- (v5);
\draw[blue, very thick] (v1) arc (90:210:.7);

\draw (-1.5,0) circle (.7cm)
	node at (-1.3,-1.2) {(b) \ \small $w\sim v_4; \ H_1$};
\coordinate (o2) at (-1.5,0);
\coordinate (a1) at (-1.5,.7);
\draw[fill=black] (a1) circle (0.15em) node at (-1.5cm,.45cm) {\small $v_1$};
\coordinate (a2) at ($(o2)!1!-60:(a1)$);
\draw[fill=black] (a2) circle (0.15em) node at (-1.1cm,.25cm) {\small $v_2$};
\coordinate (a3) at ($(o2)!1!-120:(a1)$);
\draw[fill=black] (a3) circle (0.15em);
\coordinate (a4) at ($(o2)!1!180:(a1)$);
\draw[fill=black] (a4) circle (0.15em) node at (-1.5cm,-.45cm) {\small $v_4$};
\coordinate (a5) at ($(o2)!1!120:(a1)$);
\draw[fill=black] (a5) circle (0.15em);
\coordinate (a6) at ($(o2)!1!60:(a1)$);
\draw[fill=black] (a6) circle (0.15em);
\coordinate (A) at (-.4,.7);
\draw[fill=black] (A) circle (0.15em) node[right] {\small $u$};
\coordinate (B) at (-.4,-.7);
\draw[fill=black] (B) circle (0.15em) node[right] {\small $w$};
\draw[blue, very thick] (A) -- (B);
\draw[blue, very thick] (A) -- (a1);
\draw[blue, very thick] (B) -- (a2);
\draw[blue, very thick] (B) -- (a4);
\draw[blue, very thick] (a1) arc (90:-90:.7);

\draw (1.5,0) circle (.7cm)
	node at (1.5,-1.2) {(c) \ \small $w\sim v_6; \ H_2$};
\coordinate (o3) at (1.5,0);
\coordinate (b1) at (1.5,.7);
\draw[fill=black] (b1) circle (0.15em) node at (1.5cm,.45cm) {\small $v_1$};
\coordinate (b2) at ($(o3)!1!-60:(b1)$);
\draw[fill=black] (b2) circle (0.15em) node at (1.9cm,.25cm) {\small $v_2$};
\coordinate (b3) at ($(o3)!1!-120:(b1)$);
\draw[fill=black] (b3) circle (0.15em);
\coordinate (b4) at ($(o3)!1!180:(b1)$);
\draw[fill=black] (b4) circle (0.15em);
\coordinate (b5) at ($(o3)!1!120:(b1)$);
\draw[fill=black] (b5) circle (0.15em);
\coordinate (b6) at ($(o3)!1!60:(b1)$);
\draw[fill=black] (b6) circle (0.15em) node at (1.1cm,.25cm) {\small $v_6$};
\coordinate (A) at (1.5,1.25);
\draw[fill=black] (A) circle (0.15em) node at (1.65cm,1.1cm) {\small $u$};
\coordinate (B) at (1.5,1.8);
\draw[fill=black] (B) circle (0.15em) node[right] {\small $w$};
\draw[blue, very thick] (A) -- (B);
\draw[blue, very thick] (A) -- (b1);
\draw[blue, very thick] (B) -- (b2);
\draw[blue, very thick] (B) -- (b6);
\draw[blue, very thick] (b6) arc (150:30:.7);

\draw (4.5,0) circle (.7cm)
	node at (4.6,-1.2) {(d) \ \small $|N(u)\cap A| > 2; \ H_2$};
\coordinate (o4) at (4.5,0);
\coordinate (c1) at (4.5,.7);
\draw[fill=black] (c1) circle (0.15em) node at (4.5cm,.45cm) {\small $v_1$};
\coordinate (c2) at ($(o4)!1!-60:(c1)$);
\draw[fill=black] (c2) circle (0.15em) node at (4.9cm,.25cm) {\small $v_2$};
\coordinate (c3) at ($(o4)!1!-120:(c1)$);
\draw[fill=black] (c3) circle (0.15em);
\coordinate (c4) at ($(o4)!1!180:(c1)$);
\draw[fill=black] (c4) circle (0.15em);
\coordinate (c5) at ($(o4)!1!120:(c1)$);
\draw[fill=black] (c5) circle (0.15em);
\coordinate (c6) at ($(o4)!1!60:(c1)$);
\draw[fill=black] (c6) circle (0.15em);
\coordinate (A) at ($(o4)!2!(c1)$);
\draw[fill=black] (A) circle (0.15em) node[left] {\small $u$};
\coordinate (B) at ($(o4)!2!-60:(c1)$);
\draw[fill=black] (B) circle (0.15em) node[below] {\small $w$};
\coordinate (C) at ($(c2)!.5!(A)$);
\draw[fill=black] (C) circle (0.15em) node[right] {\small $x$};
\draw[blue, very thick] (A) -- (B);
\draw[blue, very thick] (A) -- (C);
\draw[blue, very thick] (A) -- (c1);
\draw[blue, very thick] (B) -- (c2);
\draw[blue, very thick] (C) -- (c2);
\draw[blue, very thick] (c1) arc (90:30:.7);
\end{tikzpicture}
\end{subfigure}
\caption{Invalid structures in $G$ if Claim 4.10 is false.}
\end{figure}
\noindent Then for all $k\geq3$ we have that $A=A_{1,0}\cup A_{2,0}\cup A_{1,1}$, and vertices in $A_{1,1}$ which are adjacent to each other have adjacent attaching vertices in $H$. When possible, we will define $\eta(v)=\{3\}$ for vertices $v\in A$. However, this not possible on all vertices in $A_{1,1}$, nor for vertices in $A_{2,0}$ which have a neighbor in $H$ that is already assigned $\{3\}$. So define $$\overline{A}=\left\{u\in A_{2,0}:\{\tau(v):v\in N(u)\cap H\}=\{\sigma(u),\{3\}\}\right\}\cup\left\{u\in A_{1,1}:\sigma(u)=\{2\}\right\}$$
\vskip.2cm
\noindent Now define $\eta:G\rightarrow K_3$ by:\\
\\
\hspace*{1.7cm} $\eta(u)=\left\{\begin{array}{ll}
\tau(u) & \text{ if } u\in H\\
\{1\} & \text{ if } u\in \overline{A}\cap A_{1,1};\\
& \hspace{.4cm} u\in \overline{A}\cap A_{2,0} \text{ with } \sigma(u)=\{2\}\\
\{2\} & \text{ if } u\in \overline{A}\cap A_{2,0} \text{ with } \sigma(u)=\{1\}\\
\{3\} & \text{ if } u\in A\setminus\overline{A};\\ 
 & \hspace{.4cm} u\in G\setminus (A\cup H) \text{ with } N(u)\cap \overline{A}\neq\varnothing\\
\sigma(v) & \text{else}\\
 \end{array}\right\}$
\vskip.6cm
\noindent Notice that for $u\in\overline{A}$, $\eta(u)=\{1,2\}\setminus\sigma(u)=\sigma(w)$ for any $w\in N(u)$. Hence the last qualifier assigns $\{3\}$ to any neighbors of vertices in $\overline{A}$ which need to be switched, and $\eta=\sigma$ on all remaining vertices.\\
\\
It is clear that this defines a homomorphism on $A\cup H$, and on any edges between vertices $u \in G\setminus(A\cup H)$ for which $N(u)\cap \overline{A}=\varnothing$. Thus it suffices to check that this is consistent on edges with a vertex $u\in G\setminus (A\cup H)$ for which $N(u)\cap\overline{A}\neq\varnothing$. Note that by definition, any such $u$ will have $\eta(u)=\{3\}$. So to ensure that $\eta(w)\neq\{3\}$ for any $w\in N(u)$, we must check that $N(u)\cap\left(A\setminus\overline{A}\right)=\varnothing$, and for all $w\in N(u)$ that $N(w)\cap\overline{A}=\varnothing$. The cases $k=3$ and $k\geq4$ must again be handled separately.
\vskip.3cm
\begin{claim}
Let $k\geq4$, and let $u\in G\setminus (A\cup H)$ with $x\in N(u)\cap\overline{A}$. Then $N(u)\cap\left(A\setminus\overline{A}\right)=\varnothing$, and if $w\in N(u)$ then $N(w)\cap\overline{A}=\varnothing$.
\end{claim}
\vskip.2cm
\noindent {\bf Proof of Claim.} Let $u$ and $x$ be as in the claim. Since $x\in A\subseteq B_1$, it must be that $u\in (B_1\cup B_2)\setminus A$. Suppose first that $u\in B_1$. Then $|N(u)\cap B_1|\leq1$ by Claim 4.9, and hence $N(u)\cap A=\{ x\}$. So $N(u)\cap\left(A\setminus\overline{A}\right)=\varnothing$. Similarly, Claim 4.9 implies that $N(x)\cap B_1=\{ u\}$ and that $|N(x)\cap H|=1$, so $x\in A_{1,0}$, which contradicts our assumption that $x\in\overline{A}\subseteq A_{2,0}\cup A_{1,1}$. So $u\notin B_1$.\\
\\
Suppose instead that $u\in B_2$, and thus $N(u)\subseteq B_1\cup B_2\cup B_3$. For $w\in N(u)\cap B_3$ it is immediate that $w\notin A$ and that $\left(N(w)\cap\overline{A}\right)\subseteq (N(w)\cap B_1)=\varnothing$. Then assume that there is $w\in N(u)\cap B_1, \ w\neq x$. Then $\sigma(w)=\sigma(x)$ since $G$ is bipartite, and thus for any attaching points in $H$, $v_i\sim x, \ v_j\sim w$, it must be that $i$ and $j$ are either both odd or both even. Since $x\in A_{2,0}\cup A_{1,1}$ we may assume that either $i\neq j$ or $x\in A_{1,1}$. Recall that an adjacent pair in $A_{1,1}$ forms a square with its adjacent pair in $H$. So if $i=j$ and $x\in A_{1,1}$ with $\{z\}=N(x)\cap A_{1,1}$, then $z$ is adjacent to one of the two vertices in $N(v_i)\cap H$. As shown in Figure 4.6 (a), this creates a copy of $H_1$. So we may assume that $i\neq j$. But then the short path between $v_i$ and $v_j$ in $H$ along with $\{ u,w,x\}$ creates a cycle of length $l, \ 6\leq l\leq k+4$.
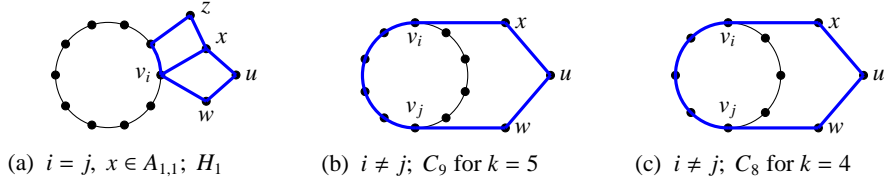
\begin{figure}[H]
\hspace{.35cm}
\begin{subfigure}[t]{.3\textwidth}
\begin{tikzpicture}
\draw (0,0) circle (.7cm)
	node at (.1,-1.2) {\small (a) \ $i=j, \ x\in A_{1,1}; \ H_1$};
\coordinate (v1) at (72:.7);
\draw[fill=black] (v1) circle (0.15em);
\coordinate (v2) at (36:.7);
\draw[fill=black] (v2) circle (0.15em);
\coordinate (v3) at (0:.7);
\draw[fill=black] (v3) circle (0.15em) node[left] {\small $v_i$};
\coordinate (v4) at (324:.7);
\draw[fill=black] (v4) circle (0.15em);
\coordinate (v5) at (288:.7);
\draw[fill=black] (v5) circle (0.15em);
\coordinate (v6) at (252:.7);
\draw[fill=black] (v6) circle (0.15em);
\coordinate (v7) at (216:.7);
\draw[fill=black] (v7) circle (0.15em);
\coordinate (v8) at (180:.7);
\draw[fill=black] (v8) circle (0.15em);
\coordinate (v9) at (144:.7);
\draw[fill=black] (v9) circle (0.15em);
\coordinate (v10) at (108:.7);
\draw[fill=black] (v10) circle (0.15em);
\coordinate (A) at (0:1.7);
\draw[fill=black] (A) circle (0.15em) node[right] {\small $u$};
\coordinate (B) at (15:1.35);
\draw[fill=black] (B) circle (0.15em) node at (1.5,.5) {\small $x$};
\coordinate (C) at (-15:1.35);
\draw[fill=black] (C) circle (0.15em) node[below] {\small $w$};
\coordinate (D) at (36:1.35);
\draw[fill=black] (D) circle (0.15em) node at (1.3,.9) {\small $z$};
\draw[blue, very thick] (A) -- (B);
\draw[blue, very thick] (A) -- (C);
\draw[blue, very thick] (B) -- (v3);
\draw[blue, very thick] (B) -- (D);
\draw[blue, very thick] (C) -- (v3);
\draw[blue, very thick] (D) -- (v2);
\draw[blue, very thick] (v2) arc (36:0:.7);
\end{tikzpicture}
\end{subfigure}
\hspace{.2cm}
\begin{subfigure}[t]{.3\textwidth}
\begin{tikzpicture}
\draw (0,0) circle (.7cm)
	node at (.2,-1.2) {\small (b) \ $i\neq j; \ C_{9}$ for $k=5$};
\coordinate (v1) at (90:.7);
\draw[fill=black] (v1) circle (0.15em) node[below] {\small $v_i$};
\coordinate (v2) at (54:.7);
\draw[fill=black] (v2) circle (0.15em);
\coordinate (v3) at (18:.7);
\draw[fill=black] (v3) circle (0.15em);
\coordinate (v4) at (342:.7);
\draw[fill=black] (v4) circle (0.15em);
\coordinate (v5) at (306:.7);
\draw[fill=black] (v5) circle (0.15em);
\coordinate (v6) at (270:.7);
\draw[fill=black] (v6) circle (0.15em) node[above] {\small $v_j$};
\coordinate (v7) at (234:.7);
\draw[fill=black] (v7) circle (0.15em);
\coordinate (v8) at (198:.7);
\draw[fill=black] (v8) circle (0.15em);
\coordinate (v9) at (162:.7);
\draw[fill=black] (v9) circle (0.15em);
\coordinate (v10) at (126:.7);
\draw[fill=black] (v10) circle (0.15em);
\coordinate (A) at (1.8,0);
\draw[fill=black] (A) circle (0.15em) node[right] {\small $u$};
\coordinate (B) at (1.2,.7);
\draw[fill=black] (B) circle (0.15em) node[right] {\small $x$};
\coordinate (C) at (1.2,-.7);
\draw[fill=black] (C) circle (0.15em) node[right] {\small $w$};
\draw[blue, very thick] (A) -- (B);
\draw[blue, very thick] (A) -- (C);
\draw[blue, very thick] (B) -- (v1);
\draw[blue, very thick] (C) -- (v6);
\draw[blue, very thick] (v1) arc (90:270:.7);
\end{tikzpicture}
\end{subfigure}
\hspace{.2cm}
\begin{subfigure}[t]{.3\textwidth}
\begin{tikzpicture}
\draw (0,0) circle (.7cm)
	node at (.2,-1.2) {\small (c) \ $i\neq j; \ C_8$ for $k=4$};
\coordinate (v1) at (90:.7);
\draw[fill=black] (v1) circle (0.15em) node[below] {\small $v_i$};
\coordinate (v2) at (45:.7);
\draw[fill=black] (v2) circle (0.15em);
\coordinate (v3) at (0:.7);
\draw[fill=black] (v3) circle (0.15em);
\coordinate (v4) at (315:.7);
\draw[fill=black] (v4) circle (0.15em);
\coordinate (v5) at (270:.7);
\draw[fill=black] (v5) circle (0.15em) node[above] {\small $v_j$};
\coordinate (v6) at (225:.7);
\draw[fill=black] (v6) circle (0.15em);
\coordinate (v7) at (180:.7);
\draw[fill=black] (v7) circle (0.15em);
\coordinate (v8) at (135:.7);
\draw[fill=black] (v8) circle (0.15em);
\coordinate (A) at (1.8,0);
\draw[fill=black] (A) circle (0.15em) node[right] {\small $u$};
\coordinate (B) at (1.2,.7);
\draw[fill=black] (B) circle (0.15em) node[right] {\small $x$};
\coordinate (C) at (1.2,-.7);
\draw[fill=black] (C) circle (0.15em) node[right] {\small $w$};
\draw[blue, very thick] (A) -- (B);
\draw[blue, very thick] (A) -- (C);
\draw[blue, very thick] (B) -- (v1);
\draw[blue, very thick] (C) -- (v5);
\draw[blue, very thick] (v1) arc (90:270:.7);
\end{tikzpicture}
\end{subfigure}
\caption{$w\in N(u)\cap B_1, \ w\neq x$.}
\end{figure}
\noindent If $k\geq 5$ this length is less than $2k$, as in Figure 4.6 (b), which contradicts the minimality of $H$. If $k=4$, this creates a cycle of length 8 only if $j=(i+4)$ mod 8, as in Figure 4.6 (c). Recall that as defined for $k=4$ (see Figure 4.2 (b)), $\tau(v_i)\neq\tau(v_{(i+4)\text{ mod } 8})$ for all $i$. But $\sigma(x)=\sigma(w)$, and therefore $x$ and $w$ cannot both be in $A$. As this is true for any pair of vertices in $N(u)\cap B_1$, it must be that $|N(u)\cap A|\leq1$, and hence by assumption $N(u)\cap A=\{ x\}$. Furthermore, for any $w\in B_1\setminus A$, if $w\sim z$ with $z\in A$, then by Claim 4.9 it must be that $z\in A_{1,0}$. Thus $N(w)\cap\overline{A}=\varnothing$ for all $w\in N(u)\cap B_1$.\\
\\
Now suppose that $w\in N(u)\cap B_2$. Note that $w\notin A\subseteq B_1$. Then let $x=x_u$ and let $x_w\in N(w)\cap B_1$. Note that $\sigma(x_u)\neq\sigma(x_w)$ since there is a path of length 3 from $x_u$ to $x_w$, so the connecting vertices for $x_u$ and $x_w$ in $H$ must be distinct to avoid creating an odd cycle. Then the path between their connecting vertices in $H$ along with $x_u,u,w$, and $x_w$ creates a cycle of length $l, \ 6\leq l\leq k+5$. For $k > 5$ such a cycle contradicts the minimality of $H$.\\
\\
Suppose $k=4$. Then to avoid creating a cycle of length less than 8, this construction requires that the path in $H$ between the connecting points of $x_u$ and $x_w$ is length 3. Suppose by way of contradiction that $x_w\in\overline{A}$. So both $x_u,x_w\in\overline{A}\subseteq A_{2,0}\cup A_{1,1}$. Without loss of generality, assume that $x_u\in A_{2,0}$ with $N(x_u)\cap H=\{ v_i, v_{(i+2)\text{ mod }8}\}$ for some $1\leq i\leq8$. Then $v_{(i+5)\text{ mod } 8}$ is the only vertex in $H$ which is distance 3 from both $v_i$ and $v_{(i+2)\text{ mod }8}$. So $x_w$ cannot also be in $A_{2,0}$, and thus $x_w\in A_{1,1}$ with $N(w_x)\cap H=\{ v_{(i+5)\text{ mod } 8}\}$. By definition, $\sigma(v)=\{2\}$ for all $v\in\overline{A}\cap A_{1,1}$, so we have that $\sigma(x_w)=\{2\}$, and thus also $\tau(v_{(i+5)\text{ mod }8})=\{2\}$. For $\tau$ as defined, this forces $i\in\{3,5,8\}$. But now $\sigma(x_u)=\{1\}$, and so $x_u\in\overline{A}\cap A_{2,0}$ implies that $\{\tau(v_i),\tau(v_{(i+2)\text{ mod }8})\}=\{\{1\},\{3\}\}$, which is only possible if $i\in\{1,4\}$. Figure 4.7 (a) and (b) shows the arrangements if $i=1$ and $i=4$ respectively, with $\tau(v_6),\tau(v_1)\neq\{2\}$. Hence it must be that both $x_u,x_w\in A_{1,1}$. However, $\sigma(x_u)\neq\sigma(x_w)$, and as noted above $\sigma(v)=\{2\}$ for all $v\in\overline{A}\cap A_{1,1}$. So at most one of $x_u$ and $x_w$ can be in $\overline{A}$, and by assumption $x_u\in\overline{A}$. Thus $x_w\notin\overline{A}$ and $N(w)\cap\overline{A}=\varnothing$ for all $w\in N(u)\cap B_2$.\\
\begin{figure}[H]
\begin{subfigure}[t]{.3\textwidth}
\begin{tikzpicture}
\draw (0,0) circle (1cm) node at (0,-2.75) {\small (a) \ $i=1, \ \tau(v_6)\neq\{2\}$};
\coordinate (v3) at (0:1);
\draw[fill=black] (v3) circle (0.15em) node[left] {\small $v_3$} node[right] {\footnotesize $\{3\}$};
\coordinate (v2) at (45:1);
\draw[fill=black] (v2) circle (0.15em) node[below left] {\small $v_2$} node[above right] {\footnotesize $\{2\}$};
\coordinate (v1) at (90:1);
\draw[fill=black] (v1) circle (0.15em) node at (0,.7) {\small $v_1$} node[above] {\footnotesize $\{1\}$};
\coordinate (v8) at (135:1);
\draw[fill=black] (v8) circle (0.15em) node[below right] {\small $v_8$} node[above left] {\footnotesize $\{2\}$};
\coordinate (v7) at (180:1);
\draw[fill=black] (v7) circle (0.15em) node[right] {\small $v_7$} node[left] {\footnotesize $\{1\}$};
\coordinate (v6) at (225:1);
\draw[fill=black] (v6) circle (0.15em) node[above right] {\small $v_6$} node[below left] {\footnotesize $\{3\}$};
\coordinate (v5) at (270:1);
\draw[fill=black] (v5) circle (0.15em) node at (0,-.7) {\small $v_5$} node[below] {\footnotesize $\{2\}$};
\coordinate (v4) at (315:1);
\draw[fill=black] (v4) circle (0.15em) node[above left] {\small $v_4$} node[below right] {\footnotesize $\{1\}$};
\coordinate (vu) at (45:2);
\draw[fill=black] (vu) circle (0.15em) node[above] {\small $x_u$};
\coordinate (vw) at (180:2);
\draw[fill=black] (vw) circle (0.15em) node[below] {\small $x_w$};
\coordinate (u) at (90:2);
\draw[fill=black] (u) circle (0.15em) node[above] {\small $u$};
\coordinate (w) at (135:2);
\draw[fill=black] (w) circle (0.15em) node[left] {\small $w$};
\draw[blue, very thick] (vu) -- (v1);
\draw[thin] (vu) -- (v3);
\draw[blue, very thick] (vw) -- (v6);
\draw[blue, very thick] (u) -- (w);
\draw[blue, very thick] (w) -- (vw);
\draw[blue, very thick] (vu) -- (u);
\draw[blue, very thick] (v1) arc (90:225:1);
\end{tikzpicture}
\end{subfigure}
\hspace{1.5cm}
\begin{subfigure}[t]{.3\textwidth}
\begin{tikzpicture}
\draw (0,0) circle (1cm) node at (0,-2.75) {\small (b) \ $i=4, \ \tau(v_1)\neq\{2\}$};
\coordinate (v3) at (0:1);
\draw[fill=black] (v3) circle (0.15em) node[left] {\small $v_3$} node[right] {\footnotesize $\{3\}$};
\coordinate (v2) at (45:1);
\draw[fill=black] (v2) circle (0.15em) node[below left] {\small $v_2$} node[above right] {\footnotesize $\{2\}$};
\coordinate (v1) at (90:1);
\draw[fill=black] (v1) circle (0.15em) node at (0,.7) {\small $v_1$} node[above] {\footnotesize $\{1\}$};
\coordinate (v8) at (135:1);
\draw[fill=black] (v8) circle (0.15em) node[below right] {\small $v_8$} node[above left] {\footnotesize $\{2\}$};
\coordinate (v7) at (180:1);
\draw[fill=black] (v7) circle (0.15em) node[right] {\small $v_7$} node[left] {\footnotesize $\{1\}$};
\coordinate (v6) at (225:1);
\draw[fill=black] (v6) circle (0.15em) node[above right] {\small $v_6$} node[below left] {\footnotesize $\{3\}$};
\coordinate (v5) at (270:1);
\draw[fill=black] (v5) circle (0.15em) node at (0,-.7) {\small $v_5$} node[below] {\footnotesize $\{2\}$};
\coordinate (v4) at (315:1);
\draw[fill=black] (v4) circle (0.15em) node[above left] {\small $v_4$} node[below right] {\footnotesize $\{1\}$};
\coordinate (vu) at (270:2);
\draw[fill=black] (vu) circle (0.15em) node[left] {\small $x_u$};
\coordinate (vw) at (45:2);
\draw[fill=black] (vw) circle (0.15em) node[right] {\small $x_w$};
\coordinate (u) at (315:2);
\draw[fill=black] (u) circle (0.15em) node[right] {\small $u$};
\coordinate (w) at (0:2);
\draw[fill=black] (w) circle (0.15em) node[right] {\small $w$};
\draw[blue, very thick] (vu) -- (v4);
\draw[thin] (vu) -- (v6);
\draw[blue, very thick] (vw) -- (v1);
\draw[blue, very thick] (u) -- (w);
\draw[blue, very thick] (w) -- (vw);
\draw[blue, very thick] (vu) -- (u);
\draw[blue, very thick] (v1) arc (90:-45:1);
\end{tikzpicture}
\end{subfigure}
\caption{$x_u\in\overline{A}\cap A_{2,0}, x_w\in\overline{A}\cap A_{1,1}, \ \sigma(x_u)=\{1\}$ when $k=4$.}
\end{figure}

\noindent Now suppose $k=5$, and again assume that $x_w\in\overline{A}$. Here cycles created by adjoining $u,w,x_u,$ and $x_w$ will be of length 10 only if the connecting vertices in $H=C_{10}$ are distance 5 apart. Hence neither $x_u$ nor $x_w$ can be in $A_{2,0}$, and both must be in $A_{1,1}$. Thus, as in the case above for $k=4$, $\sigma(x_u)\neq\sigma(x_w)$ implies that at most one of $x_u$ and $x_w$ can be in $\overline{A}$. Therefore $x_w\notin\overline{A}$ and $N(w)\cap\overline{A}=\varnothing$ for all $w\in N(u)\cap B_2$. $_{\square}$
\vskip.4cm
\begin{claim}
Let $k=3$, and let $u\in G\setminus (A\cup H)$ with $x\in N(u)\cap\overline{A}$. Then $N(u)\cap\left(A\setminus\overline{A}\right)=\varnothing$, and if $w\in N(u)$ then $N(w)\cap\overline{A}=\varnothing$.
\end{claim}
\vskip.3cm
\noindent {\bf Proof of Claim.} Let vertices $u$ and $x$ be as in the claim, and let $v_i\in N(x)\cap H$ be a vertex such that $\sigma(x)=\tau(v_i)$, which exists by the assumption that $x\in A$. Then since $G$ is bipartite, it must be that $\{v\in H: d(v,u)=2\}\subseteq\{ v_i,v_{(i+2)\text{ mod }6},v_{(i+4)\text{ mod }6}\}$. Recall that $\tau$ is defined such that for each $j\in K_3$ the set $\tau^{-1}(\{ j\})$ contains one vertex with an odd index and one vertex with an even index (see Figure 4.2 (a)). So $\{ v_i,v_{(i+2)\text{ mod }6},v_{(i+4)\text{ mod }6}\}\cap\tau^{-1}(\sigma(x))=\{ v_i\}$.  Hence $w\sim v_i$ for every $w\in N(u)\cap A$. Then if $|N(u)\cap A|\geq3$, this would create a copy of $H_2$, depicted in Figure 4.8 (a). Thus $|N(u)\cap A|\leq 2$.\\
\begin{figure}[H]
\hspace*{.7cm}
\begin{subfigure}[t]{\textwidth}
\begin{tikzpicture}
\draw (-4,0) circle (.7cm)
	node at (-4,-1.2) {\footnotesize (a) \ $|N(u)\cap A|\geq 3; \ H_2$};
\coordinate (o1) at (-4,0);
\coordinate (v1) at (-4,.7);
\draw[fill=black] (v1) circle (0.15em);
\coordinate (v2) at ($(o1)!1!-60:(v1)$);
\draw[fill=black] (v2) circle (0.15em) node at (-3.65,.25) {\small $v_i$};
\coordinate (v3) at ($(o1)!1!-120:(v1)$);
\draw[fill=black] (v3) circle (0.15em);
\coordinate (v4) at ($(o1)!1!180:(v1)$);
\draw[fill=black] (v4) circle (0.15em);
\coordinate (v5) at ($(o1)!1!120:(v1)$);
\draw[fill=black] (v5) circle (0.15em);
\coordinate (v6) at ($(o1)!1!60:(v1)$);
\draw[fill=black] (v6) circle (0.15em);
\coordinate (A) at (-2.3,.1);
\draw[fill=black] (A) circle (0.15em) node[right] {\small $u$};
\coordinate (B) at (-2.8,.7);
\draw[fill=black] (B) circle (0.15em) node[right] {\small $x$};
\coordinate (C) at (-2.8,.2);
\draw[fill=black] (C) circle (0.15em);
\coordinate (D) at (-2.8,-.3);
\draw[fill=black] (D) circle (0.15em);
\draw[blue, very thick] (A) -- (B);
\draw[blue, very thick] (A) -- (C);
\draw[blue, very thick] (A) -- (D);
\draw[blue, very thick] (B) -- (v2);
\draw[blue, very thick] (C) -- (v2);
\draw[blue, very thick] (D) -- (v2);

\draw (0,0) circle (.7cm)
	node at (-.2,-1.2) {\footnotesize (b) \ $|N(u)\cap A|=2$,};
\coordinate (o2) at (0,0) node at (.4,-1.6) {\footnotesize $|N(u)\cap A_{2,0}|\geq1; \ H_1$};
\coordinate (a1) at (0,.7);
\draw[fill=black] (a1) circle (0.15em);
\coordinate (a2) at ($(o2)!1!-60:(a1)$);
\draw[fill=black] (a2) circle (0.15em) node at (.35,.25) {\small $v_i$};
\coordinate (a3) at ($(o2)!1!-120:(a1)$);
\draw[fill=black] (a3) circle (0.15em);
\coordinate (a4) at ($(o2)!1!180:(a1)$);
\draw[fill=black] (a4) circle (0.15em);
\coordinate (a5) at ($(o2)!1!120:(a1)$);
\draw[fill=black] (a5) circle (0.15em);
\coordinate (a6) at ($(o2)!1!60:(a1)$);
\draw[fill=black] (a6) circle (0.15em);
\coordinate (A) at (1.7,0);
\draw[fill=black] (A) circle (0.15em) node[right] {\small $u$};
\coordinate (B) at (1.2,.7);
\draw[fill=black] (B) circle (0.15em) node[right] {\small $x$};
\coordinate (C) at (1.2,-.7);
\draw[fill=black] (C) circle (0.15em);
\draw[blue, very thick] (A) -- (B);
\draw[blue, very thick] (A) -- (C);
\draw[blue, very thick] (B) -- (a2);
\draw[blue, very thick] (C) -- (a2);
\draw[blue, very thick] (C) -- (a4);
\draw[blue, very thick] (a4) arc (-90:30:.7);

\draw (4,0) circle (.7cm)
	node at (3.8,-1.2) {\footnotesize (c) \ $|N(u)\cap A|=2,$};
\coordinate (o3) at (4,0) node at (4.4,-1.6) {\footnotesize $|N(u)\cap A_{1,1}|\geq1; \ H_1$};
\coordinate (b1) at (4,.7);
\draw[fill=black] (b1) circle (0.15em);
\coordinate (b2) at ($(o3)!1!-60:(b1)$);
\draw[fill=black] (b2) circle (0.15em) node at (4.35,.25) {\small $v_i$};
\coordinate (b3) at ($(o3)!1!-120:(b1)$);
\draw[fill=black] (b3) circle (0.15em);
\coordinate (b4) at ($(o3)!1!180:(b1)$);
\draw[fill=black] (b4) circle (0.15em);
\coordinate (b5) at ($(o3)!1!120:(b1)$);
\draw[fill=black] (b5) circle (0.15em);
\coordinate (b6) at ($(o3)!1!60:(b1)$);
\draw[fill=black] (b6) circle (0.15em);
\coordinate (A) at (5.7,.5);
\draw[fill=black] (A) circle (0.15em) node[right] {\small $u$};
\coordinate (B) at (5.2,.7);
\draw[fill=black] (B) circle (0.15em) node at (5.4,.8) {\small $x$};
\coordinate (C) at (5.2,.2);
\draw[fill=black] (C) circle (0.15em) node at (5.4,.1) {\small $w$};
\coordinate (D) at (5.2,-.4);
\draw[fill=black] (D) circle (0.15em) node[right] {\small $w'$};
\draw[blue, very thick] (A) -- (B);
\draw[blue, very thick] (A) -- (C);
\draw[blue, very thick] (B) -- (b2);
\draw[blue, very thick] (C) -- (b2);
\draw[blue, very thick] (C) -- (D);
\draw[blue, very thick] (D) -- (b3);
\draw[blue, very thick] (b3) arc (-30:30:.7);
\end{tikzpicture}
\end{subfigure}
\caption{$k=3 \text{ and } |N(u)\cap A|\geq2$.}
\end{figure}

\noindent Suppose that $|N(u)\cap A|=2$. If there was an adjacency from the set $\{ v_{(i+2)\text{ mod }6},v_{(i+4)\text{ mod }6}\}$ to $N(u)\cap A$, it would create a copy of $H_1$, as in Figure 4.8 (b). So $N(u)\cap A_{2,0}=\varnothing$. And adjacent pairs in $A_{1,1}$ form squares with their connecting pair in $H$. So if $w\in N(u)\cap A_{1,1}$, then there is $w'\in N(w)\cap A_{1,1}$ such that, without loss of generality, $w'\sim v_{i+1 \text{ mod } 6}$. But, as depicted in Figure 4.8 (c), this creates a copy of $H_1$. Thus it cannot be that that both $|N(u)\cap A| > 1$ and $N(u)\cap (A_{2,0}\cup A_{1,1})\neq\varnothing$. Since $x\in\overline{A}\subseteq A_{2,0}\cup A_{1,1}$ by assumption, we have that $N(u)\cap A=\{ x\}$, and so $N(u)\cap\left(A\setminus\overline{A}\right)=\varnothing$.\\
\\
Now it suffices to check for all $w\in N(u)\setminus (H\cup\{ x\})$ \ that \ $N(w)\cap\overline{A}=\varnothing$. Suppose that $u\in B_1$. Then $N(u)\cap H\subset\{ v_{(i+1)\text{ mod }6}, v_{(i+3)\text{ mod }6}, v_{(i+5)\text{ mod }6}\}$ since $x\sim v_i$. Further suppose that $x\in A_{2,0}$, so without loss of generality $N(x)\cap H=\{ v_i, v_{(i+2)\text{ mod }6}\}$. Then all three possible adjacencies for $u$ in $H$ are shown in Figure 4.9, with $u\sim v_{(i+1)\text{ mod }6}$ creating a copy of $H_2$, while $v\sim v_{(i+3)\text{ mod }6}$ \ or \ $v\sim v_{(i+5)\text{ mod }6}$ each creates a copy of $H_1$. Hence $x\notin A_{2,0}$.
\begin{figure}[H]
\hspace*{.7cm}
\begin{subfigure}[t]{.3\textwidth}
\begin{tikzpicture}
\draw (0,0) circle (.7cm)
	node at (0,-1.2) {\small (a) \ $u\sim v_{(i+1)\text{ mod }6}; \ H_2$};
\coordinate (v1) at (90:.7);
\draw[fill=black] (v1) circle (0.15em) node at (0,.42) {\small $v_i$};
\coordinate (v2) at (30:.7);
\draw[fill=black] (v2) circle (0.15em);
\coordinate (v3) at (330:.7);
\draw[fill=black] (v3) circle (0.15em);
\coordinate (v4) at (270:.7);
\draw[fill=black] (v4) circle (0.15em);
\coordinate (v5) at (210:.7);
\draw[fill=black] (v5) circle (0.15em);
\coordinate (v6) at (150:.7);
\draw[fill=black] (v6) circle (0.15em);
\coordinate (A) at (1,.5);
\draw[fill=black] (A) circle (0.15em) node at (.9,.3) {\small $u$};
\coordinate (B) at (1.5,.7);
\draw[fill=black] (B) circle (0.15em) node[right] {\small $x$};
\draw[blue, very thick] (A) -- (B);
\draw[blue, very thick] (A) -- (v2);
\draw[blue, very thick] (B) -- (v1);
\draw[blue, very thick] (B) -- (v3);
\draw[blue, very thick] (v1) arc (90:-30:.7);
\end{tikzpicture}
\end{subfigure}
\hspace{.2cm}
\begin{subfigure}[t]{.3\textwidth}
\begin{tikzpicture}
\draw (0,0) circle (.7cm)
	node at (0.3,-1.2) {\small (b) \ $u\sim v_{(i+3)\text{ mod }6}; \ H_1$};
\coordinate (v1) at (90:.7);
\draw[fill=black] (v1) circle (0.15em) node at (0,.42) {\small $v_i$};
\coordinate (v2) at (30:.7);
\draw[fill=black] (v2) circle (0.15em);
\coordinate (v3) at (330:.7);
\draw[fill=black] (v3) circle (0.15em);
\coordinate (v4) at (270:.7);
\draw[fill=black] (v4) circle (0.15em);
\coordinate (v5) at (210:.7);
\draw[fill=black] (v5) circle (0.15em);
\coordinate (v6) at (150:.7);
\draw[fill=black] (v6) circle (0.15em);
\coordinate (A) at (1.3,-.7);
\draw[fill=black] (A) circle (0.15em) node[right] {\small $u$};
\coordinate (B) at (1.3,.7);
\draw[fill=black] (B) circle (0.15em) node[right] {\small $x$};
\draw[blue, very thick] (A) -- (B);
\draw[blue, very thick] (A) -- (v4);
\draw[blue, very thick] (B) -- (v1);
\draw[blue, very thick] (B) -- (v3);
\draw[blue, very thick] (v1) arc (90:-90:.7);
\end{tikzpicture}
\end{subfigure}
\begin{subfigure}[t]{.3\textwidth}
\begin{tikzpicture}
\draw (0,0) circle (.7cm)
	node at (0.3,-1.2) {\small (c) \ $u\sim v_{(i+5)\text{ mod }6}; \ H_1$};
\coordinate (v1) at (90:.7);
\draw[fill=black] (v1) circle (0.15em);
\coordinate (v2) at (30:.7);
\draw[fill=black] (v2) circle (0.15em) node at (.35,.25) {\small $v_i$};
\coordinate (v3) at (330:.7);
\draw[fill=black] (v3) circle (0.15em);
\coordinate (v4) at (270:.7);
\draw[fill=black] (v4) circle (0.15em);
\coordinate (v5) at (210:.7);
\draw[fill=black] (v5) circle (0.15em);
\coordinate (v6) at (150:.7);
\draw[fill=black] (v6) circle (0.15em);
\coordinate (A) at (1.3,.7);
\draw[fill=black] (A) circle (0.15em) node[right] {\small $u$};
\coordinate (B) at (1.3,-.7);
\draw[fill=black] (B) circle (0.15em) node[right] {\small $x$};
\draw[blue, very thick] (A) -- (B);
\draw[blue, very thick] (A) -- (v1);
\draw[blue, very thick] (B) -- (v2);
\draw[blue, very thick] (B) -- (v4);
\draw[blue, very thick] (v1) arc (90:-90:.7);
\end{tikzpicture}
\end{subfigure}
\caption{$k=3, \ u\in B_1, \text{ and } x\in A_{2,0}$.}
\end{figure}
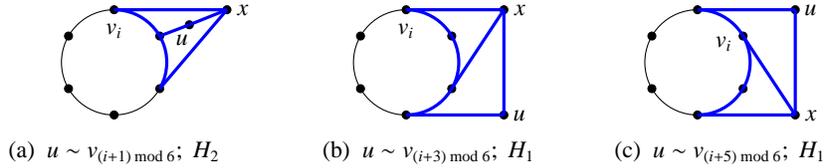

\noindent So $x\in\overline{A}\cap A_{1,1}$ which implies that $\sigma(x)=\{2\}$, and so $\sigma(u)=\{1\}$. Then without loss of generality let $i=2$. Note $x\sim v_2$ implies $u$ must be connected to a vertex in $H$ which has odd index. But $u\notin A$, so $u\nsim v_1$. Let $\{ w\}=N(x)\cap A_{1,1}$, so it must be that $\sigma(w)=\{1\}$ and $w\sim v_1$. Then $u\sim v_3$ would create a copy of $H_1$, shown in Figure 4.10 (a), and so $u\sim v_5$. Now consider a vertex $z\in N(u)\setminus(H\cup\{ x\})$, and suppose there exists $v_z\in N(z)\cap\overline{A}$. Then $d(u,v_z)=2$, so $\sigma(v_z)=\sigma(u)=\{1\}$. And by definition, we have that $\{ v\in\overline{A}:\sigma(v)=\{1\}\}\subseteq A_{2,0}$, so $v_z\in A_{2,0}$ and the two vertices in $N(v_z)\cap H$ have images $\{1\}$ and $\{3\}$ under $\sigma$. Since $\sigma(v_z)=\sigma(u)$, vertices in $N(v_z)\cap H$ must have odd index. Thus $N(v_z)\cap H=\{ v_1,v_3\}$. As shown in Figure 4.10 (b), $\{ w,x,v_z,v_1,v_2,v_3\}$ is then a copy of $H_1$. So if $u\in B_1$ there cannot be any vertex $z\in N(u)$ with $N(z)\cap\overline{A}\neq\varnothing$.\\
\begin{figure}[H]
\hspace{.1cm}
\begin{subfigure}[t]{.3\textwidth}
\begin{tikzpicture}
\draw (0,0) circle (1cm)
	node at (0,-1.5) {\small (a) \ $u\sim v_3; \ H_1$};
\coordinate (v6) at (150:1);
\draw[fill=black] (v6) circle (0.15em) node at (-.6,.4) {\small $v_6$};
\coordinate (v1) at (90:1);
\draw[fill=black] (v1) circle (0.15em) node at (0,.75) {\small $v_1$};
\coordinate (v2) at (30:1);
\draw[fill=black] (v2) circle (0.15em) node at (.6,.4) {\small $v_2$};
\coordinate (v3) at (330:1);
\draw[fill=black] (v3) circle (0.15em) node at (.6,-.4) {\small $v_3$};
\coordinate (v4) at (270:1);
\draw[fill=black] (v4) circle (0.15em) node at (0,-.75) {\small $v_4$};
\coordinate (v5) at (210:1);
\draw[fill=black] (v5) circle (0.15em) node at (-.6,-.4) {\small $v_5$};
\coordinate (O) at (0,0);
\coordinate (A) at ($(O)!1.75!(v1)$);
\draw[fill=black] (A) circle (0.15em) node[left] {\small $w$};
\coordinate (B) at ($(O)!1.75!(v2)$);
\draw[fill=black] (B) circle (0.15em) node[right] {\small $x$};
\coordinate (C) at ($(O)!1.75!(v3)$);
\draw[fill=black] (C) circle (0.15em) node[right] {\small $u$};
\draw[blue, very thick] (A) -- (B);
\draw[blue, very thick] (A) -- (v1);
\draw[blue, very thick] (B) -- (v2);
\draw[blue, very thick] (B) -- (C);
\draw[blue, very thick] (C) -- (v3);
\draw[blue, very thick] (v1) arc (90:-30:1);
\end{tikzpicture}
\end{subfigure}
\hspace{.8cm}
\begin{subfigure}[t]{.3\textwidth}
\begin{tikzpicture}
\draw (0,0) circle (1cm)
	node at (0,-1.5) {\small (b) \ $u\sim z\sim v_z, \ v_z\in\overline{A}; \ H_1$};
\coordinate (v6) at (150:1);
\draw[fill=black] (v6) circle (0.15em) node at (-.6,.4) {\small $v_6$};
\coordinate (v1) at (90:1);
\draw[fill=black] (v1) circle (0.15em) node at (0,.75) {\small $v_1$};
\coordinate (v2) at (30:1);
\draw[fill=black] (v2) circle (0.15em) node at (1.1,.55) {\small $v_2$};
\coordinate (v3) at (330:1);
\draw[fill=black] (v3) circle (0.15em) node at (.6,-.4) {\small $v_3$};
\coordinate (v4) at (270:1);
\draw[fill=black] (v4) circle (0.15em) node at (0,-.75) {\small $v_4$};
\coordinate (v5) at (210:1);
\draw[fill=black] (v5) circle (0.15em) node at (-.6,-.4) {\small $v_5$};
\coordinate (A) at (0,1.4);
\draw[fill=black] (A) circle (0.15em) node[left] {\small $w$};
\coordinate (B) at ($(v3)!1.87!(v2)$);
\draw[fill=black] (B) circle (0.15em) node[right] {\small $x$};
\coordinate (C) at (-2,1.8);
\draw[fill=black] (C) circle (0.15em) node[left] {\small $u$};
\coordinate (E) at (1.5,-.5);
\draw[fill=black] (E) circle (0.15em) node[right] {\small $v_z$};
\coordinate (D) at (1.5,1.7);
\draw[fill=black] (D) circle (0.15em) node[right] {\small $z$};
\draw[blue, very thick] (A) -- (B);
\draw[blue, very thick] (A) -- (v1);
\draw[blue, very thick] (B) -- (v2);
\draw[black, thin] (B) -- (C);
\draw[black, thin] (C) -- (D);
\draw[black, thin] (D) -- (E);
\draw[blue, very thick] (E) -- (v1);
\draw[blue, very thick] (E) -- (v3);
\draw[blue, very thick] (v1) arc (90:-30:1);
\draw[thin] (C) -- (v5);
\end{tikzpicture}
\end{subfigure}
\caption{$k=3, \ u\in B_1, \text{ and } x\in\overline{A}\cap A_{1,1}$.}
\end{figure}

\noindent Now assume that $u\in B_2$. Let $x=x_u$, and suppose that there is a vertex $w\in N(u)$ which has a vertex $x_w\in N(w)\cap\overline{A}$. Then there is a path of length 3 from $x_u$ to $x_w$, so $x_u\neq x_w$ and also $\sigma(x_u)\neq\sigma(x_w)$. Then they cannot both be in $A_{1,1}$, since $\sigma(v)=\{2\}$ for all $v\in\overline{A}\cap A_{1,1}$. Without loss of generality, let $x_u\in A_{2,0}$. Suppose that $x_w\in A_{1,1}$ and let $\{ z\}=N(x_w)\cap A_{1,1}$. Then $\sigma(x_w)=\{2\}$ and without loss of generality we have that $x_w\sim v_2$ and $z\sim v_1$. Now $x_u$ must be adjacent to two vertices in $H$ with odd indices, one of which has image $\{3\}$ under $\tau$. And $\sigma(x_u)=\{1\}$, so $N(x_u)\cap H=\{ v_1,v_3\}$. This is depicted in Figure 4.11 (a), where $\{ v_1,v_2,v_3,x_u,x_w,z\}$ forms a copy of $H_1$. Suppose instead that $x_w\in A_{2,0}$. Then each $x_w$ and $x_u$ must be adjacent to a vertex whose image under $\tau$ is $\{3\}$. Without loss of generality assume that $x_u\sim v_3$ and $x_w\sim v_6$, and let $x_u\sim v_j$ and $x_w\sim v_l$ be the other two adjacencies. Then $j$ must be odd and $l$ must be even, and $\tau(v_j)=\sigma(x_u)\neq\sigma(x_w)=\tau(v_l)$. So either $j=1$ and $l=2$, or $j=5$ and $l=4$. By symmetry we may assume that $x_u\sim v_1$ and $x_w\sim v_2$, as in Figure 4.11 (b). But then $\{ v_1,v_2,v_3,v_6,x_u,x_w\}$ forms a copy of $H_1$. Thus $N(w)\cap\overline{A}=\varnothing$ for all $w\in N(u)$. $_{\square}$
\begin{figure}[H]
\hspace*{1cm}
\begin{subfigure}[t]{.35\textwidth}
\begin{tikzpicture}
\draw (0,0) circle (1cm)
	node at (0,-1.5) {\small (a) \ $x_w\in A_{1,1}; \ H_1$};
\coordinate (v6) at (150:1);
\draw[fill=black] (v6) circle (0.15em) node at (-.6,.4) {\small $v_6$};
\coordinate (v1) at (90:1);
\draw[fill=black] (v1) circle (0.15em) node at (0,.75) {\small $v_1$};
\coordinate (v2) at (30:1);
\draw[fill=black] (v2) circle (0.15em) node at (1.1,.55) {\small $v_2$};
\coordinate (v3) at (330:1);
\draw[fill=black] (v3) circle (0.15em) node at (.6,-.4) {\small $v_3$};
\coordinate (v4) at (270:1);
\draw[fill=black] (v4) circle (0.15em) node at (0,-.75) {\small $v_4$};
\coordinate (v5) at (210:1);
\draw[fill=black] (v5) circle (0.15em) node at (-.6,-.4) {\small $v_5$};
\coordinate (O) at (0,0);
\coordinate (A) at (0,1.9);
\draw[fill=black] (A) circle (0.15em) node[left] {\small $z$};
\coordinate (B) at (.86,1.9);
\draw[fill=black] (B) circle (0.15em) node[above] {\small $x_w$};
\coordinate (C) at (1.7,1.9);
\draw[fill=black] (C) circle (0.15em) node[right] {\small $w$};
\coordinate (D) at (1.7,.5);
\draw[fill=black] (D) circle (0.15em) node[right] {\small $u$};
\coordinate (E) at (1.7,-.5);
\draw[fill=black] (E) circle (0.15em) node[right] {\small $x_u$};
\draw[blue, very thick] (A) -- (B);
\draw[blue, very thick] (A) -- (v1);
\draw[blue, very thick] (B) -- (v2);
\draw[thin] (B) -- (C);
\draw[thin] (C) -- (D);
\draw[thin] (D) -- (E);
\draw[blue, very thick] (E) -- (v1);
\draw[blue, very thick] (E) -- (v3);
\draw[blue, very thick] (v1) arc (90:-30:1);
\end{tikzpicture}
\end{subfigure}
\hspace{1cm}
\begin{subfigure}[t]{.35\textwidth}
\begin{tikzpicture}
\draw (0,0) circle (1cm)
	node at (0,-1.5) {\small (b) \ $x_w\in A_{2,0}; \ H_1$};
\coordinate (v6) at (150:1);
\draw[fill=black] (v6) circle (0.15em) node at (-.6,.4) {\small $v_6$};
\coordinate (v1) at (90:1);
\draw[fill=black] (v1) circle (0.15em) node at (0,.75) {\small $v_1$};
\coordinate (v2) at (30:1);
\draw[fill=black] (v2) circle (0.15em) node at (1.1,.55) {\small $v_2$};
\coordinate (v3) at (330:1);
\draw[fill=black] (v3) circle (0.15em) node at (.6,-.4) {\small $v_3$};
\coordinate (v4) at (270:1);
\draw[fill=black] (v4) circle (0.15em) node at (0,-.75) {\small $v_4$};
\coordinate (v5) at (210:1);
\draw[fill=black] (v5) circle (0.15em) node at (-.6,-.4) {\small $v_5$};
\coordinate (O) at (0,0);
\coordinate (A) at (0,1.9);
\draw[fill=black] (A) circle (0.15em) node[above] {\small $x_w$};
\coordinate (B) at (1.7,1.9);
\draw[fill=black] (B) circle (0.15em) node[right] {\small $w$};
\coordinate (C) at (1.7,.5);
\draw[fill=black] (C) circle (0.15em) node[right] {\small $u$};
\coordinate (D) at (1.7,-.5);
\draw[fill=black] (D) circle (0.15em) node[right] {\small $x_u$};
\draw[blue, very thick] (A) -- (v2);
\draw[blue, very thick] (A) -- (v6);
\draw[thin] (A) -- (B);
\draw[thin] (B) -- (C);
\draw[thin] (C) -- (D);
\draw[blue, very thick] (D) -- (v1);
\draw[blue, very thick] (D) -- (v3);
\draw[blue, very thick] (v6) arc (150:-30:1);
\end{tikzpicture}
\end{subfigure}
\caption{$k=3, \ u\in B_2, \text{ and } \{x_u,x_w\}\cap A_{2,0}\neq\varnothing$.}
\end{figure}
\vskip.3cm
\noindent It follows immediately from Claims 4.11 and 4.12 that $\eta:G\rightarrow K_3$ is a homomorphism for all $k\geq3$. And $\eta$ extends $\tau$, so $\eta_{\{ l,j\}}$ extends $\tau_{\{ l,j\}}$. Thus by Lemma 4.2, $\eta$ and $\eta_{\{ l,j\}}$ are in distinct connected components of Hom$(G,K_3)$, which is therefore disconnected. $_{\square}$
\vskip.5cm
The details of this proof suggests that it is only necessary to require that there be some cycle $C_{2k}\subset G, \ k\geq3$, such that $H'\cap C_{2k}=\varnothing$ for any subgraph $H'\subset G$ with $H'\in\{ H_1,H_2\}$. But even this leaves out many bipartite $G$ for which Hom$(G,K_3)$ is disconnected. For example, consider the graph $G'$ which is a circular ladder with six rungs. So $G'$ is two disjoint copies of $C_6$ with respective vertices labeled consecutively as $v_1,\ldots, v_6$, and $w_1,\ldots w_6$, and additional edges $v_i\sim w_i$ for each $i$. Note that $G'$ is bipartite, does not admit any folds, and every edge is contained in at least two copies of $H_1$. Define $\eta(v_i)=\{ i\text{ mod }3\}$ and $\eta(w_i)=\{(i+1)\text{ mod }3\}$ for all $i$. Then let $\tau$ be the restriction of $\eta$ to the induced cycle on $\{ v_1,\ldots,v_6\}$. So $r(\tau)=0$, and as before $\eta$ and $\eta_{\{ l,j\}}$ are in distinct connected components of Hom$(G',K_3)$ for any fixed pair $\{ l,j\}\subset\{1,2,3\}$.
\par In the contrarian case of $Q_3$, again every edge is contained in multiple copies of $H_1$. The real issue, however, is that every 6-cycle in $Q_3$ is contained in either an induced copy of $H_1$ or an induced copy of $Q_3\setminus\{ v\}$ for some $v\in Q_3$. Both of these graphs fold to an edge, so Hom$(H_1,K_3)$ and Hom$(Q_3\setminus\{ v\},K_3)$ are each connected. But then any distinct 0-cells $\tau_1,\tau_2\in$ Hom$(C_6,K_3)$ which extend to all of $Q_3$ must first extend to either $H_1$ or $Q_3\setminus\{ v\}$, where their extensions are necessarily in a single connected component.

\section{Hom$\left(G(n,p),K_m\right)$}

The topological connectivity of random hom-complexes has been studied previously under the guise of the Neighborhood Complex of $G(n,p)$, the Erd\H{o}s--R\'{e}nyi model for a random graph. $G(n,p)$ is the probability space of graphs on $n$ vertices where each edge is inserted independently with probability $p=p(n)$. Kahle ~\cite{Kahle1} showed that the connectivity of $\mathcal{N}(G(n,p))$ is concentrated between $1/2$ and $2/3$ of the expected value of the largest clique in $G(n,p)$, and also obtained asymptotic bounds on the number of dimensions with non-trivial homology. As we discussed earlier, $\mathcal{N}(G(n,p))$ is homotopy equivalent to Hom$(K_2,G(n,p))$. Here we take the opposite perspective and consider the random polyhedral complex Hom$(G(n,p),K_m)$.
\par The major benefit of utilizing $D(G)$ in a generalization of the \u{C}uki\'{c}-Kozlov Theorem is that $k$-cores have been well-studied in random graphs models. We say that $G(n,p)$ has property $\mathcal{P}$ {\it with high probability} if $\displaystyle\lim_{n\rightarrow\infty}\text{Pr}[G(n,p)\in\mathcal{P}]=1$. A property $\mathcal{P}$ is said to have a sharp threshold $\hat{p}=\hat{p}(n)$ if for all $\epsilon > 0$ $$\lim_{n\rightarrow\infty}\text{Pr}[G(n,p)\in\mathcal{P}]=\left\{\begin{array}{ll}
0 & \text{ if } p\leq(1-\epsilon)\hat{p}\\
1 & \text{ if } p\geq(1+\epsilon)\hat{p}\end{array}\right\}$$
\par When $p=c/n$ for constant $c > 0$, Pittel, Spencer and Wormwald ~\cite{PSWcore} showed that the existence of a $k$-core in $G(n,p)$ has a sharp threshold $c=c_k$ for $k\geq3$, and that asymptotically $c_k=k+\sqrt{k\log k}+O(\log k)$. Approximate values for small $k$ are known, such as $c_3\approx3.35, \ c_4\approx5.14, \ c_5\approx6.81$. When $k=2$, the existence of cycles has a one-sided sharp threshold at $c=1$. For $c > 1$, indeed $\text{Pr}[D(G(n,c/n))\geq2]\rightarrow1$. However, for all $0 < c < 1$ there is a constant $0 < f(c) < 1$ for which $\text{Pr}[D(G(n,c/n))\geq2]\rightarrow f(c)$. To simplify this issue, define $c_k$ for all $k\geq2$ by $$c_k\coloneqq\sup\left\{c > 0:\lim_{n\rightarrow\infty}\text{Pr}[D(G(n,c/n)\geq k]=0\right\}$$ In particular $c_2=0$, and for $k\geq 3$ these are precisely the sharp thresholds mentioned above.
\par By applying Theorem 2.8 to these thresholds, we then immediately obtain lower bounds on the topological connectivity of Hom$(G(n,p),K_m)$. For notational convenience, define $M(n,c,m)\coloneqq conn\left[\text{Hom}(G(n,c/n),K_m)\right]$.
\begin{thm}
If $k\geq2$ and $p=c/n$ with $c < c_{k+1}$, then for all $m\geq3$ $$\lim_{n\rightarrow\infty}\text{Pr}\left[M(n,c,m)\geq m-k-2\right]=1$$
\end{thm}
Theorem 2.8 does not require that the input graph be connected, and in fact $G(n,p)$ will be disconnected when $p=c/n$. Let the disjoint connected components of $G(n,p)$ be $G_0,G_1,\ldots,G_t$, ordered from most vertices to least. For $c > 1$, with high probability $G_0$ is a giant component containing more than half the vertices, and $G_i$ is either an isolated vertex, a tree, or a unicyclic graph for $i\geq1$. So $$\text{Hom}(G(n,c/n),K_m)=\text{Hom}\left(\coprod_0^{t}G_i,K_m\right)=\prod_{i=0}^{t}\text{Hom}(G_i,K_m)$$ And thus $$M(n,c,m)=\min_{0\leq i\leq t}\left\{ conn[{\text{Hom}(G_i,K_m)}]\right\}$$ 
\par For $G_i=\{v_i\}$, Hom$(G_i,K_m)=\Delta^2$, which is contractible. If a finite connected graph $G_i$ is a tree, then it folds to a single edge and Hom$(G_i,K_m)\simeq\text{Hom}(K_2,K_m)\simeq S^{m-2}$. Finally, if $G_i$ is a finite unicyclic graph, then it folds to its cycle $C_n$. If $n=4$, then $C_4$ folds further to $K_2$ and Hom$(G_i,K_m)\simeq S^{m-2}$. Thus, with high probability we have $$\text{Hom}(G(n,c/n),K_m)=\left(\Delta^2\right)^{t_1}\times\left(S^{m-2}\right)^{t_2}\times\text{Hom}(G_0,K_m)\times\prod_{j=1}^{t_3}\text{Hom}(C_{n_j},K_m)$$ where $t_1=$ the number of isolated vertices, $t_2=$ the number of trees and unicyclic components whose cycle is $C_4$, $t_3=t-(t_1+t_2)$, and $n_j\neq 4$ for all $j$. 
\par When $n_j=3$, Hom$(C_{n_j},K_m)=\text{Hom}(K_3,K_m)$, which is a wedge of $(m-3)$-spheres ~\cite{BabKoz}. Kozlov also computed the homology of Hom$(C_n,K_m)$ for all $n\geq 5, \ m\geq 4$ in ~\cite{Koz1}. In particular, for all $m\geq 4$:\\
\hspace*{3.2cm} $conn\left[\text{Hom}\left(C_{2r+1},K_m\right)\right]=m-4$ \ \ for $r\geq1$\\
\hspace*{3.2cm} $conn\left[\text{Hom}\left(C_{2r},K_m\right)\right]=m-3$ \ \ for $r\geq2$\\
\\
Note that the cases $n=3$ and $n=4$, which yield spheres, are consistent with these values. And the threshold for all small components to be isolated vertices is $p=\frac{\log n}{4n}$, so for $c < c_{k+1}, \ k\geq2$, with high probability there is some $i\geq1$ such that $$conn[\text{Hom}(G_i,K_m)]\in\{ m-3,m-4\}$$
This provides an upper bound for $M(n,c,m)$, and so when $m\geq4$ we improve Theorem 5.1 to the following:
\begin{thm}
If $k\geq2$ and $p=c/n$ with $c < c_{k+1}$, then for all $m\geq 4$ $$\lim_{n\rightarrow\infty}\text{Pr}\left[m-k-2\leq M(n,c,m)\leq m-3 \ \right]=1$$
\end{thm}
\vskip.2cm
For fixed $m$, the lower bound decreases as $c$ gets bigger, and the gap between upper and lower bounds becomes worse. We expect that $conn\left[\text{Hom}(G_0,K_m)\right]$ should decrease as $D(G_0)$ increases and the giant component becomes more highly connected. However, the lower bound in Theorem 2.8 is not tight in general. For example, the complete bipartite graph $K_{i,j}$ folds to an edge, so Hom$(K_{i,j},K_m)\simeq S^{m-2}$. But $D(K_{i,j})=\min\{i,j\}$, so Theorem 2.8 yields $conn\left[\text{Hom}(K_{i,j},K_m)\right]\geq m-\min\{i,j\}-2$. If $i,j\geq2$, then this bound is not tight, and for large $i,j$ it can be arbitrarily bad.
\par To sharpen these results, we turn to examining the chromatic number of $G(n,c/n)$. Similar to the thresholds for the appearance of $k$-cores, Achlioptas and Friedgut ~\cite{AchFried1} showed that there is a sharp threshold sequence $d_k(n)$ such that for any $\epsilon > 0$, $$c < (1-\epsilon)d_k(n) \ \Longrightarrow \lim_{n\rightarrow\infty}\text{Pr}[\chiup(G(n,c/n))\leq k]=1$$ $$c > (1+\epsilon)d_k(n) \ \Longrightarrow\lim_{n\rightarrow\infty}\text{Pr}[\chiup(G(n,c/n))\geq k+1]=1$$ The convergence of the $d_k(n)$ remains an open problem, but Coja-Oghlan and Vilenchik ~\cite{C-OV1} improved the lower bound on $\displaystyle\liminf_{n\rightarrow\infty}d_k(n)$ to within a small constant of the upper bound for $\displaystyle\limsup_{n\rightarrow\infty}d_k(n)$. In particular, for $o_k(1)$ a term which goes to 0 as $k$ becomes large, $$d_{k}^{-}\coloneqq\liminf_{n\rightarrow\infty}d_k(n)\geq 2k\log k-\log k-2\log2-o_k(1)$$ $$d_{k}^{+}\coloneqq\limsup_{n\rightarrow\infty}d_k(n)\leq2k\log k-\log k-1+o_k(1)$$ Then as $k$ increases, the difference $|d_{k}^{+}-d_{k}^{-}|$ approaches $2\log2-1\approx0.39$. And in the context of hom-complexes it is immediate that with high probability $$\text{Hom}(G(n,c/n),K_m)=\varnothing \ \text{ for } c > d_{m}^{+} \ \text{ and } \ \text{Hom}(G(n,c/n),K_m)\neq\varnothing \ \text{ for } c < d_{m}^{-}$$
For $m=3$, we can then evaluate the connectivity of Hom$(G(n,c/n),K_3)$ for all $c > 0$, excluding the gap between $d_3^{-}$ and $d_3^{+}$.
\vskip.3cm
\begin{thm}
For $\displaystyle 1\leq c < d_3^{-}, \ \text{Hom}(G(n,c/n),K_3)$ is disconnected with high probability.
\end{thm}
\vskip.2cm
\noindent {\bf Proof of Theorem 5.3.} By the definition of $d_3^{-}$ and the existence of odd cycles with high probability for $c\geq1$, $\chiup(G(n,c/n))=3$. Then by Theorem 4.6 $\text{Hom}(G(n,c/n),K_3)$ is disconnected with high probability. $_{\square}$
\vskip.3cm
When $0 < c < 1$, all connected components are isolated vertices, trees, or unicyclic graphs with high probability. If there is a component which contains edges and does not fold to a single edge, then Hom$(G(n,c/n)$ will be disconnected. But if every connected component of $G(n,c/n)$ is an isolated vertex or folds to an edge, then $\text{Hom}(G(n,c/n),K_3)$ is connected. In the latter case, specifically $$\text{Hom}(G(n,c/n),K_3)=\left(\Delta^2\right)^{t_1}\times\mathbb{T}^{t_2}$$ where $t_1=$ number of isolated vertices, and $t_2=$ the number of connected components which contain at least one edge. And $t_2 > 0$ with high probability for $p=c/n$, so the complex will not be contractible, and $M(n,c,3)=0$. 
\begin{thm}
For $\displaystyle 0 < c < 1, \ c'=\frac{1}{2}\log(1-c)+\frac{c}{2}+\frac{c^2}{4}+\frac{c^4}{8}$, $$\lim_{n\rightarrow\infty}\text{Pr}[M(n,c,3)=-1]=1-e^{c'} \ \text{ and } \ \lim_{n\rightarrow\infty}\text{Pr}[M(n,c,3)=0]=e^{c'}$$
\end{thm}
\noindent {\bf Proof of Theorem 5.4.} If $G$ is a connected unicyclic graph whose cycle is $C_l$ for $l\neq4$, then $G$ folds to $C_l$ and Hom$(G,K_3)\simeq\text{ Hom}(C_l,K_3)$, which is disconnected. So if $G(n,c/n)$ does not contain any cycle $C_l$ for $l\neq4$, then all components of $G(n,c/n)$ will be isolated vertices or will fold to an edge. For a fixed $l$, the number of $l$-cycles in $G(n,c/n)$ approaches a limiting Poisson distribution with mean $\frac{c^l}{2l}$, and so $$\text{Pr}[C_l\nsubset G(n,c/n) \ \text{ for all } l\neq4]\rightarrow\exp\left\{-\frac{c^3}{6}-\sum_{l=5}^{\infty}\frac{c^l}{2l}\right\}$$ For $c < 1$, the sum in the exponent converges to $\displaystyle c'=\frac{1}{2}\log(1-c)+\frac{c}{2}+\frac{c^2}{4}+\frac{c^4}{8}$. Thus $\text{Pr}[M(c,n,3)=0]\rightarrow e^{c'}$ and $\text{Pr}[M(c,n,3)=-1]\rightarrow1-e^{c'}$. $_{\square}$\\
\vskip.5cm
\section{Further Questions}
The biggest question left unanswered is that of bounding $conn\left[\text{Hom}(G_0,K_m)\right]$ from above. We expect that $conn\left[\text{Hom}(G_0,K_m)\right]$ remains close to $m-k-2$ when this bound makes sense, that $M(n,c,m)=conn\left[\text{Hom}(G_0,K_m)\right]$, and that $M(n,c,m)$ is a non-increasing function for fixed $m$, with high probability. However, exhibiting non-trivial homology classes in Hom$(G,H)$ is difficult in general. In Section 4 we showed that disconnected components can be lifted via subgraphs, which can be viewed as lifting non-trivial 0-cycles, but there is no known analogous result for lifting higher dimensional non-trivial cycles.
\par We also expect that $M(n,c,m)$ increases monotonically when $c$ is fixed and $m$ is increasing, but this is not true in general for a fixed input graph. For instance, consider again our favorite counter example, $Q_3$. We noted in Section 4 that $conn\left[\text{Hom}(Q_3,K_3)\right]\geq0$. And by Theorem 2.8, $conn\left[\text{Hom}(Q_3,K_5)\right]\geq 0$. But for $m=4$, there are colorings of $Q_3$ in which each color class is a pair of antipodal corners. In such a coloring, every vertex is adjacent to one vertex from each of the other three color class, so it represents an isolated 0-cell in $\text{Hom}(Q_3,K_4)$. Thus $conn\left[\text{Hom}(Q_3,K_4)\right]=-1$. In this example with $Q_3$, the problem seems to occur when $m\leq D(G)$, and it is plausible that $conn\left[\text{Hom}(G,K_m)\right]$ is either increasing or non-decreasing for $m\geq D(G)+1$, if $G$ is a fixed graph which cannot be reduced via folds. This, however, remains an open question.
\par In the random setting, considering $m\leq D(G)$ presents a serious roadblock to evaluating $M(n,c,m)$ for all pairs of $c$ and $m$. Since the $d_k^{\pm}$ grow much faster than the $c_k$, there are intervals where $D(G(n,c/n))$ is much larger than $\chiup(G(n,c/n))$. For example, Molloy ~\cite{Mol1} pointed out that $c_{40}\approx52.23$ while $d_{19} > 53.88$. So if $c=53$, with high probability $D(G(n,c/n))=40$, but $m$ can be as small as 19 before Hom$(G(n,c/n),K_m)$ becomes empty. When $m\leq D(G)$ the lower bound from Theorem 2.8 provides no information, and new methods are required.
\par A different direction altogether would be to indulge in a closer examination of precise numerical estimates on the Betti numbers and Euler characteristic of Hom$(G(n,c/n),K_m)$. \u{C}uki\'{c} and Kozlov's ~\cite{CukKoz2} work on cycles in hom-complexes makes the case that $m=3$ a tantalizingly tractable place to start this type of investigation.
\\

\bibliographystyle{plain}
\bibliography{Master}

\begin{thebibliography}{10}

\bibitem{AchFried1}
Dimitris Achlioptas and Ehud Friedgut.
\newblock A sharp threshold for {$k$}-colorability.
\newblock {\em Random Structures Algorithms}, 14(1):63--70, 1999.

\bibitem{BabKoz}
Eric Babson and Dmitry~N. Kozlov.
\newblock Complexes of graph homomorphisms.
\newblock {\em Israel J. Math.}, 152:285--312, 2006.

\bibitem{CervdHJ}
Luis Cereceda, Jan van~den Heuvel, and Matthew Johnson.
\newblock Mixing 3-colourings in bipartite graphs.
\newblock {\em European J. Combin.}, 30(7):1593--1606, 2009.

\bibitem{C-OV1}
Amin Coja-Oghlan and Dan Vilenchik.
\newblock Chasing the {$k$}-colorability threshold.
\newblock In {\em 2013 {IEEE} 54th {A}nnual {S}ymposium on {F}oundations of
  {C}omputer {S}cience---{FOCS} 2013}, pages 380--389. IEEE Computer Soc., Los
  Alamitos, CA, 2013.

\bibitem{Csorba1}
Peter Csorba.
\newblock {\em Non-Tidy Spaces and Graph Colorings}.
\newblock PhD thesis, ETH Z\"{u}rich, 2005.

\bibitem{CukKoz1}
Sonja~Lj. {\v{C}}uki{\'c} and Dmitry~N. Kozlov.
\newblock Higher connectivity of graph coloring complexes.
\newblock {\em Int. Math. Res. Not.}, (25):1543--1562, 2005.

\bibitem{CukKoz2}
Sonja~Lj. {\v{C}}uki{\'c} and Dmitry~N. Kozlov.
\newblock The homotopy type of complexes of graph homomorphisms between cycles.
\newblock {\em Discrete Comput. Geom.}, 36(2):313--329, 2006.

\bibitem{Eng1}
Alexander Engstr{\"o}m.
\newblock A short proof of a conjecture on the connectivity of graph coloring
  complexes.
\newblock {\em Proc. Amer. Math. Soc.}, 134(12):3703--3705 (electronic), 2006.

\bibitem{Kahle1}
Matthew Kahle.
\newblock The neighborhood complex of a random graph.
\newblock {\em J. Combin. Theory Ser. A}, 114(2):380--387, 2007.

\bibitem{Koz1}
Dmitry~N. Kozlov.
\newblock Cohomology of colorings of cycles.
\newblock {\em Amer. J. Math.}, 130(3):829--857, 2008.

\bibitem{Lov1}
L.~Lov{\'a}sz.
\newblock Kneser's conjecture, chromatic number, and homotopy.
\newblock {\em J. Combin. Theory Ser. A}, 25(3):319--324, 1978.

\bibitem{Mol1}
Michael Molloy.
\newblock A gap between the appearances of a {$k$}-core and a
  {$(k+1)$}-chromatic graph.
\newblock {\em Random Structures Algorithms}, 8(2):159--160, 1996.

\bibitem{PSWcore}
Boris Pittel, Joel Spencer, and Nicholas Wormald.
\newblock Sudden emergence of a giant {$k$}-core in a random graph.
\newblock {\em J. Combin. Theory Ser. B}, 67(1):111--151, 1996.

\bibitem{Schultz1}
Carsten Schultz.
\newblock Small models of graph colouring manifolds and the {S}tiefel manifolds
  {${\rm Hom}(C_5,K_n)$}.
\newblock {\em J. Combin. Theory Ser. A}, 115(1):84--104, 2008.

\end{thebibliography}
\end{document}